\documentclass[10pt,onecolumn,twoside,draftclsnofoot]{IEEEtran}

\usepackage{newtxtext}
\usepackage{amsmath,amsthm}
\usepackage[cal=cm,scr=euler,bb=boondox,frak=euler]{mathalfa}
\usepackage{mathtools}
\usepackage{amssymb}
\usepackage{bm}
\usepackage{siunitx}
\usepackage{enumitem}
\usepackage{pgfplots}
\pgfplotsset{compat=newest}
\usepackage[utf8]{inputenc}

\newcommand*{\subparagraph}{}  
\usepackage[compact]{titlesec}

\usepackage[font=small,skip=5pt,format=hang]{caption}
\usepackage[backend=bibtex8,style=ieee,url=false,doi=false,isbn=false]{biblatex}
\usetikzlibrary{shapes,arrows,positioning,backgrounds,matrix}
\usetikzlibrary{external}

\theoremstyle{plain}
\newtheorem{theorem}{Theorem}
\newtheorem{corollary}[theorem]{Corollary}
\newtheorem{lemma}[theorem]{Lemma}
\theoremstyle{definition}

\newtheorem{definition}[theorem]{Definition}

\DeclareMathOperator{\diag}{diag}
\DeclareMathOperator{\diam}{diam}

\DeclareMathOperator{\rank}{rank}

\DeclarePairedDelimiter\abs{\lvert}{\rvert}
\DeclarePairedDelimiter\norm{\lVert}{\rVert}
\DeclarePairedDelimiter\ip{\langle}{\rangle}

\newcommand*{\real}{\ensuremath{\mathbb{R}}}
\newcommand*{\cmplx}{\ensuremath{\mathbb{C}}}
\newcommand*{\intgr}{\ensuremath{\mathbb{Z}}}
\newcommand*{\nat}{\ensuremath{\mathbb{N}}}
\newcommand*{\torus}{\ensuremath{\mathbb{T}}}

\newcommand*{\ver}{\ensuremath{\mathscr{V}}}
\newcommand*{\edg}{\ensuremath{\mathscr{E}}}

\newcommand*{\E}{\ensuremath{\mathscr{S}}}
\newcommand*{\I}{\ensuremath{\mathcal{I}}}
\newcommand*{\M}{\ensuremath{\mathscr{M}}}
\newcommand*{\N}{\ensuremath{\mathscr{N}}}
\newcommand*{\G}{\ensuremath{\mathbb{G}}}
\newcommand*{\kk}{\ensuremath{\kappa}}

\newcommand*{\g}{\ensuremath{\sigma}}
\newcommand*{\h}{\ensuremath{\chi}}
\newcommand*{\rr}{\ensuremath{\zeta}}
\newcommand*{\gv}{\ensuremath{\Sigma}}
\newcommand*{\gstar}{\ensuremath{\varpi}}
\newcommand*{\e}{\ensuremath{e}}
\newcommand*{\w}{\ensuremath{e}}
\newcommand*{\ii}{\ensuremath{\mathfrak{j}}\hspace{0.1ex}}
\newcommand*{\unit}{\ensuremath{\partial}}
\newcommand*{\domg}{\ensuremath{\mathscr{U}}}
\newcommand*{\tr}{{\ensuremath{\scriptscriptstyle T}}}
\newcommand*{\inv}{{\ensuremath{\scriptscriptstyle -1}}}
\newcommand*{\1}{\ensuremath{\mathbf{1}}}
\newcommand*{\0}{\ensuremath{\mathbf{0}}}
\newcommand*{\n}{\mkern-1.5mu}

\newcommand*{\ddotvarphi}{\ensuremath{\skew{3}{\ddot}{\varphi}}}
\newcommand*{\dotgamma}{\ensuremath{\skew{4}{\dot}{\gamma}}}
\newcommand*{\dotx}{\ensuremath{\skew{3}{\dot}{x}}}
\newcommand*{\doty}{\ensuremath{\skew{3}{\dot}{y}}}

\newcommand*{\dotU}{\ensuremath{\skew{3}{\dot}{U}}}

\title{Distributed Synchronization of Heterogeneous\\
  Oscillators on Networks with Arbitrary Topology}
\author{Enrique~Mallada,\,\, Randy~A.~Freeman,\,\, and\, Ao~Tang%
    \thanks{E. Mallada is with Caltech, CMS Dept. 1200 E California Blvd,
    Pasadena, CA 91125, USA, \texttt{mallada@caltech.edu}.}
    \thanks{R.~Freeman is with Northwestern U., EECS
    Dept., 2145 Sheridan Rd., Evanston, IL 60208-3118, USA,
    \texttt{freeman@eecs.northwestern.edu}.  His work was supported in
    part by a grant from the Office of Naval Research.}
    \thanks{A. Tang is with Cornell U., School of ECE,
    Ithaca, NY 14853, USA, \texttt{atang@ece.cornell.edu}.
    His work was supported in part by a grant from the Office of
    Naval Research.}  
    }

\begin{document}


\maketitle
\begin{abstract}
  Many network applications rely on the synchronization of coupled
  oscillators.  For example, such synchronization can provide
  networked devices with a common temporal reference necessary for
  coordinating actions or decoding transmitted messages.  In this
  paper, we study the problem of using distributed control to achieve
  both phase and frequency synchronization of a network of coupled
  heterogeneous nonlinear oscillators.  Not only do our controllers
  guarantee zero phase error in steady state under arbitrary frequency
  heterogeneity, but they also require little knowledge of the
  oscillator nonlinearities and network topology.  Furthermore, we
  provide a global convergence analysis, in the absence of noise 
  and propagation delay, for the resulting nonlinear system whose 
  phase vector evolves on the $\bm{n}$-torus.
\end{abstract}


\section{Introduction}
\label{sec:introduction}

Achieving temporal coordination among different networked devices is a
fundamental requirement for the successful operation of many
engineering systems.  For example, it is necessary in communication
systems for recovering transmitted
messages~\cite{bregni1998historical}, in sensor networks for
coordinating wake up cycles~\cite{ye_medium_2004} or achieving
temporal measurement coherence~\cite{sundararaman2005clock}, and in
computer networks for preserving the causality of distributed
events~\cite{mallada_skewless_2013}.  Almost ubiquitously, such
coordination is accomplished by providing each node of the network
with its own local oscillator and then compensating its phase and
frequency (using information received from other devices on the
network) to achieve a common temporal reference.


Legacy applications such as public switched telephone networks and
cellular networks use a centralized hierarchical synchronization
scheme with high-precision oscillators having relative frequency
errors ranging from $0.01$ to $4.6$ parts per million
(ppm)~\cite{rec1998g,itu:g.811}.  For several reasons, however, these
traditional synchronization architectures have become increasingly
unsuitable for newer applications such as wireless sensor networks.
For example, traditional methods can break down with the failure of
only a few nodes.  In addition, many newer applications use
inexpensive oscillators having errors as high as
\SI{100}{ppm}~\cite{karl2007protocols}.  Thus, a synchronization
protocol designed for these newer applications should satisfy two
essential requirements: it should be distributed and independent of
the network topology (each node should use only its neighbors'
oscillator information to adjust it own oscillator), and it should be
robust to wide variations and uncertainty in the specifications of the
oscillators used throughout the network.



A variety of synchronization algorithms have been proposed along these
lines, jointly inspired by collective synchronization in physics and
biology~\cite{winfree1967biological,kuramoto1984cooperative,
mirollo1990synchronization, kuramoto2003chemical} and cooperative 
control in engineering networks~\cite{jadbabaie2003coordination, 
ren2007information}. One possible solution is to use monotonically 
increasing time sources (e.g. clocks) and update their times based on 
offset information~\cite{tong_theoretical_1998,simeone_pulse-coupled_2007,
simeone_physical-layer_2008, rentel_mutual_2008, schenato_average_2011,carchischzam08,carzam10,cardelzam11,
Carli:2014gd} to achieve a common {\it absolute} time reference (clock 
synchronization). This is suitable for applications in computer networks where 
a reference to an absolute time is needed (e.g. distributed databases). 
Another solution is to use periodic time sources (e.g. oscillators) interconnected with 
phase comparators~\cite{simspabarstr08, carareto_architectures_2012,lun11}
 or pulse-coupling~\cite{hong2003time, lucarelli2004decentralized,
 wang2011pco}, where the objective is to achieve common {\it relative} time 
 reference (phase synchronization) that allows temporal coordination within 
 the network (e.g. waking up simultaneously).\footnote{Although one can use
  clock synchronization  to achieve phase synchronization by simply 
  mapping the linear times onto the circle using a modulo operator, this 
  approach can lead to undesirable transients if the phases keep wrapping 
  around the circle as the linear times synchronize.  For this reason we 
  consider these as separate synchronization problems, each suitable for 
  different application areas.}


While the theoretical study of clock synchronization is fairly mature, with 
solutions that can provide zero offset error synchronization on networks with
arbitrary heterogeneous frequencies~\cite{mallada_skewless_2013} 
and asynchronous updates~\cite{Carli:2014gd}, little is known about the phase 
synchronization counterpart. For example, most phase synchronization 
solutions present nonzero steady state phase differences in the presence 
of frequency heterogeneity~\cite{simspabarstr08,carareto_architectures_2012,
hong2003time,lucarelli2004decentralized, wang2011pco} with convergence guarantees
limited to idealized scenarios such as homogeneous frequencies
~\cite{mallada2013synchronization}. The only exception is~\cite{lun11} which can 
 guarantee phase synchronization for complete graph topologies. Thus whether 
or not such systems can synchronize for arbitrary networks and arbitrary frequency
heterogeneity has remained as an open question~\cite{carareto_architectures_2012}.

In this paper, we provide a positive answer to this question under very general  
conditions. We propose two distributed controllers that can achieve phase synchronization for 
a network of arbitrarily interconnected oscillators, under mild assumptions on the 
oscillator and phase comparator characteristics. 
For example, we allow the instantaneous frequency of each oscillator to be a highly
uncertain nonlinear function of the local control input, a model consistent with most
analog oscillators (such as voltage-controlled oscillators or CMOS oscillators).
Also, unlike existing work, we allow the set of oscillator frequencies 
to be bounded, so that each oscillator may operate within a prescribed frequency 
range, even during the transient part of the response. 
Finally, we allow flexibility in the choice of the phase
comparator responses, rather than assuming as in ~\cite{lun11} that they are
sinusoidal. We only require that the measured phase difference is noiseless and
can be obtained without propagation delay.


The main contribution of the paper is a novel nonlinear convergence analysis that 
leverages recent results on the stability of equilibria of homogeneous-frequency 
coupled oscillators~\cite{mallada2013synchronization}. In particular, our controllers
are based on a Hamiltonian dynamic system defined on the graph in which each 
local minimum of the energy function represents a synchronized trajectory.  Each 
controller employs a different mechanism to dissipate energy and thereby converge
to a synchronized solution.  Furthermore, we show that any trajectories that are 
synchronized in frequency but not in phase must be unstable.

\section{Notation and terminology}
\label{sec:terminology-notation}


We let $\torus=\real/2\pi\intgr$ denote the unit circle, regarded as
the Lie group of angle addition.  We equip~$\torus$ with the usual
Riemannian metric which defines the distance between two points to be
the length of the shorter of the two arcs connecting them (so that
$\diam(\torus)=\pi$).  For $p\in\nat$, we let $\torus^p$ denote the
Cartesian product of~$p$ circles.  For $i\in\{1,\dotsc,p\}$, we let
$\unit_i$ denote the unit vector field pointing in the
counterclockwise direction on the $i^{\text{th}}$ factor of~$\torus^p$
(which we write simply as~$\unit$ when $p=1$).  Because these unit
vector fields form an ordered basis for the tangent space
of~$\torus^p$ at each point, we can represent tangent vectors
for~$\torus^p$ as elements of~$\real^p$, that is, as coordinate
vectors with respect to this basis.  Moreover, all Jacobian matrices
of mappings defined on~$\torus^p$ will be representations of the
differential with respect to this basis.

All graphs in this paper will be simple, undirected, connected graphs
having $n$ vertices (with $2\leqslant n<\infty$) and~$m$ edges (with
$m\geqslant n-1$).  We represent such a graph~$\G$ as a pair
$\G=(\ver,\edg)$ for a vertex set~$\ver$ and edge set~$\edg$.  We
label and order the vertices and edges, writing $\ver=\{1,\dotsc,n\}$
and $\edg=\{1,\dotsc,m\}$, where each edge $k\n\in\n\edg$ is an
unordered pair of distinct vertices $k=\{i,j\}\subset\ver$.  For each
vertex $i\n\in\n\ver$, we let $\N_i$ denote the following indexed set
of neighbors of~$i$:
\begin{align}
  \label{eq:1}
  \N_i&=\bigl\{(j,k)\in\ver\n\times\n\edg
  \,:\,k=\{i,j\}\bigr\}\,.
\end{align}
Thus $(j,k)\in\N_i$ if and only if $(i,k)\in\N_j$, that is, if and
only if edge~$k$ connects vertices~$i$ and~$j$.

\section{Problem statement and results}
\label{sec:new-stuff}

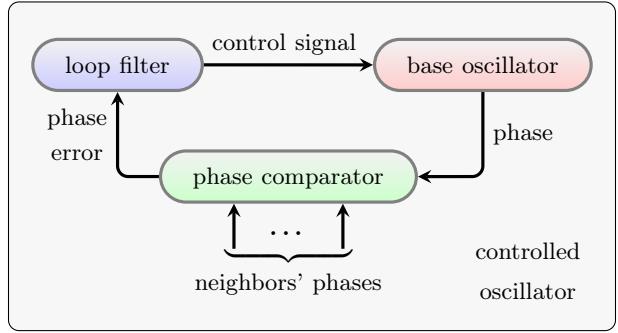
\begin{figure}
  \centering
\begin{tikzpicture}[scale=0.7,%
  component/.style={%
    text height=1.6ex,text depth=0.25ex,
    rounded rectangle,
    minimum size=6mm,inner xsep=4mm,inner ysep=2mm,
    very thick,draw=black!50,
    top color=white!95!black,bottom color=black!20},
  link/.style={->,very thick,>=stealth,rounded corners},
  align=center,
  background rectangle/.style={%
    draw=black,fill=black!3,rounded corners},
  framed,inner frame sep=2ex]
  \tikzstyle{every node}=[font=\small]
  \node (mid) {};
  \node (pc) [component,bottom color=green!20,below=of mid]
             {phase comparator};
  \node (lf) [component,bottom color=blue!20,left=of mid]
             {loop filter};
  \node (osc) [component,bottom color=red!20,right=of mid]
              {base oscillator};
  \path (lf) edge[link] node[above] {control signal\,\,} (osc);
  \draw [link] (osc) -- (osc |- pc) node[midway,right] {phase} -- (pc);
  \draw [link] (pc) -- (pc -| lf) -- (lf)
             node[midway,left] {phase\\[-0.5ex]error};
  \node (nbr1) [below left=6mm and 6mm of pc.south] {};
  \draw [link] (nbr1) -- (nbr1 |- pc.south);
  \node (nbr2) [below right=6mm and 6mm of pc.south] {};
  \draw [link] (nbr2) -- (nbr2 |- pc.south);
  \node [below=2mm of pc] {$\bm{\cdots}$};
  \node [below=2.6mm of pc,below delimiter=\},text width=10ex] {};
  \node [below=8mm of pc] {neighbors' phases};
  \node [below right=4mm and 10mm of pc]
    {\textbf{controlled}\\\textbf{oscillator}};
\end{tikzpicture}
  \caption{PLL components of a controlled oscillator.}
  \label{fig:1}
\end{figure}

We consider a network of controlled oscillators in which each
oscillator shares the current value of its phase with its immediate
neighbors.  The purpose of the controller design is to guarantee both
frequency and phase synchronization of the interconnected system.  We
adopt the classical phase-locked loop (PLL) structure for each
controlled oscillator \cite{gar79,simspabarstr08}.  This structure
consists of three components connected in feedback, as illustrated in
Fig.~\ref{fig:1}: a base oscillator, a phase comparator, and a loop
filter.  The base oscillator is a physical device (such as a
voltage-controlled oscillator) whose frequency is determined
dynamically by means of a control signal; the phase of each controlled
oscillator system is simply the phase of its base oscillator.  The
phase comparator produces a phase error signal by comparing its own
phase with the phases of its neighbors.  Finally, the loop filter
produces the control signal from the phase error.

We represent the network of oscillators by a graph $\G=(\ver,\edg)$ in
which each vertex is a controlled oscillator and each edge indicates
an exchange of phases between neighboring vertices.  We represent the
phase of oscillator~$i$ at time~$t\n\in\n\real$ by
$\varphi_i(t)\n\in\n\torus$, and we write the phase vector signal for
the entire network as the column vector
\begin{align}
  \label{eq:2}
  \varphi&=\bigl[\varphi_1\,\ldots\,\varphi_n\bigr]^\tr
  \in\,\torus^n\,.
\end{align}
We define a synchronization measure $\delta\varphi$
to be the diameter of the finite subset of~$\torus$ consisting of
the~$n$ phases~$\varphi_i$:
\begin{align}
  \label{eq:3}
  \delta\varphi
  &=\diam\bigl\{\varphi_1,\,\dotsc,\,\varphi_n\bigr\}\in[0,\pi]\,.
\end{align}
Note that if $\delta\varphi<\frac{2\pi}{3}$, then the~$n$
phases~$\varphi_i$ all lie within an arc of length $\delta\varphi$.
Therefore a small value of $\delta\varphi$ represents a tight
clustering of the oscillator phases, and these phases are all
identical when $\delta\varphi=0$.  This diameter~\eqref{eq:3} is a
worst-case measure of synchronization, rather than an average measure,
as a single outlier can make this diameter large.  The following
definition describes our design goal:
\begin{definition}\label{def:1}
  Let $\I\n\subset\n\real$ be an interval.  The network of
  oscillators achieves \textit{almost-global synchronization
    within\/}~$\I$ when for almost all initial states (including both
  the initial oscillator phases and any initial loop filter states),
  the trajectories of the system satisfy:
\begin{description}
\item[asymptotic frequency synchronization:] there exists a constant
  $\omega^\star\n\in\I$ such that
  $\dot{\varphi}_i(t)\to\omega^\star$ as $t\to\infty$ for each
  $i\n\in\n\ver$;
\item[asymptotic phase synchronization:] $\delta\varphi(t)\to 0$ as
  $t\to\infty$;
\item[constrained frequencies:] $\dot{\varphi}_i(t)\n\in\n\I$ for all
  $t\geqslant 0$ and each $i\n\in\n\ver$;
\item[internal boundedness:] all loop filter states (if any) are
  bounded in forward time.
\end{description}
\end{definition}
The constraint interval~$\I$ used in this definition characterizes the
desired range of frequencies for the oscillators.  Choosing an
appropriate interval~$\I$ in the design is thus useful for preventing
the oscillators from moving too fast, too slow, or reversing
direction.

\subsection{Base oscillator and phase comparator}
\label{sec:base-oscill-phase}

In common models of controlled analog oscillators, the instantaneous
frequency is some uncertain monotonic nonlinear function of the
control input.  For example, the graph of this function is the
``tuning curve'' seen on data sheets of many voltage-controlled
oscillators, where it is understood that this curve is typical rather
than exact.  Motivated by such models, we assign an uncertain
\textit{frequency function\/} $\h_i$ to the base oscillator in each
vertex~$i$, where $\h_i:\domg\to\real$ is a strictly increasing
$C^1$ function on a known open interval domain $\domg\n\subset\n\real$.
The angular velocity of the phase of the oscillator at time~$t$ (that
is, its instantaneous frequency) is then given by $\h_i(u_i(t))$,
where $u_i(\cdot)$ is a $\domg$-valued control signal.  Hence our
model for the base oscillator is the differential equation
\begin{align}
  \dot{\varphi}_i&=\h_i(u_i)\,.\label{eq:4}
\end{align}
In this oscillator model, the oscillator stops whenever $\h_i(u_i)=0$
and reverses direction whenever $\h_i(u_i)$ changes sign.  However, an
oscillator need not admit such behavior, as the image $\h_i(\domg)$
need not contain zero.

As we see from the model~\eqref{eq:4}, the base oscillators are
nonlinear and heterogeneous.  In addition, we do not need precise
knowledge of the frequency functions~$\h_i$ to complete our design and
guarantee almost-global synchronization.  In fact, as we will see when
we list our assumptions in Section~\ref{sec:assumpt-main-result}
below, all we need to know about the functions~$\h_i$ in the
unconstrained case $\I=\real$ is that they are~$C^1$ with positive
derivatives on~$\domg$, and that the intersection of their images
$\h_i(\domg)$ is nonempty.

The phase comparator in vertex~$i$ calculates a linear combination of
functions of phase differences to produce a phase error~$\e_i$:
\begin{align}
  \label{eq:5}
  \e_i(\varphi)&=\smashoperator{\sum_{(j,k)\in\N_i}}
   a_{k}f(\varphi_j-\varphi_i)\,,
\end{align}
where the constants $a_{k}$ are positive edge weights and
$f:\torus\to\real$ is the \textit{phase coupling function}.  For
example, the coupling function used in many classical PLL designs is
the sine function $f=\sin$.

\subsection{Example: the Kuramoto model}
\label{sec:standard-model}

Suppose that the frequency function~$\h_i$ for each base
oscillator~$i$ has the simple affine form $\h_i(u_i)=\omega_i+u_i$,
where the constant $\omega_i\n\in\n\real$ represents the nominal
oscillator frequency.  Suppose further that the phase coupling
function is the sine function $f=\sin$ and that all loop filters are
constant unity gains so that $u_i\equiv\e_i$.  Then the controlled
phase equation becomes
\begin{align}
  \label{eq:6}
  \dot{\varphi}_i&=\omega_i+\smashoperator{\sum_{(j,k)\in\N_i}}
  a_k\sin(\varphi_j-\varphi_i)\,.
\end{align}
This model of a network of coupled oscillators has been studied
extensively, and we refer the reader to the survey paper
\cite{dorbul13}.  In particular, if the graph~$\G$ is complete, and if
all edge weights~$a_k$ are the same, then this is the famous Kuramoto
model of coupled oscillators~\cite{kur75}.

The oscillator network characterized by the dynamics in~\eqref{eq:6}
fails to meet the design goal of Definition~\ref{def:1}.  Indeed, the
existence of a synchronized trajectory for the model~\eqref{eq:6} in
which $\delta\varphi\equiv 0$ implies that the nominal
frequencies~$\omega_i$ are all identical.  If these frequencies were
known precisely, then we could simply cancel them out via control by
setting $u_i=e_i+\omega^\star-\omega_i$ to obtain a model of the
form~\eqref{eq:6} in which $\omega_i\n=\n\omega^\star$ for all
$i\n\in\!\ver$.  However, even in this case of identical frequencies,
almost-global synchronization has been proved only for special classes
of connected graphs, such as complete graphs or trees~\cite{dorbul13}.
Instead, we are interested in oscillators having unknown heterogeneous
frequencies on arbitrary connected graphs.  In this case one can
choose sufficiently large edge weights~$a_k$ to guarantee
``practical'' phase synchronization in which $\delta\varphi$
becomes small, provided it does not start off too
large~\cite{dorbul13}.  However, choosing large edge weights makes it
less likely that the frequencies $\dot{\varphi}_i$ will be
constrained to a desired interval~$\I$ during the transient.  In any
case, we see that we must depart from this standard model~\eqref{eq:6}
to meet the design goal of Definition~\ref{def:1}, and we do so
by choosing a non-sinusoidal coupling function~$f$ in~\eqref{eq:5}
and a nonlinear dynamic loop filter.

\subsection{Toward synchronization: a Hamiltonian system}
\label{sec:towards-synchr-hamil}

The loop filter in vertex~$i$ produces the control signal $u_i$ from
the phase error $\e_i$.  It is well known that introducing integral
action into the loop filter can compensate offset mismatches for
networks of heterogeneous oscillators
\cite{gar79,simspabarstr08,Carli:2014gd,carchischzam08,cardelzam11,
lun11}.  Thus as
a first attempt at achieving almost-global synchronization of the
oscillator network, we simply make the loop filter a scaled
integrator:
\begin{align}
  \label{eq:7}
  \dotgamma_i&=c_i e_i(\varphi)\\
  u_i&=\rr(\gamma_i)\,,\label{eq:8}
\end{align}
where~$\gamma_i(\cdot)$ is the real-valued internal filter state,
$c_i$ is a positive parameter, and $\rr:\real\to\domg$ is a scaling
function which squeezes the value of~$\gamma_i$ into the
domain~$\domg$ of the frequency function~$\h_i$.  We do not assume
that the oscillators have access to the global time variable~$t$.  As
a result, the differential $dt$ used in the construction of
$d\gamma_i/dt$ in~\eqref{eq:7} is unknown, and we account for this by
assuming that the positive constant~$c_i$ is unknown.

The system resulting from~\eqref{eq:4} and
\eqref{eq:7}--\eqref{eq:8} is
\begin{align}
  \label{eq:9}
  \dot{\varphi}_i&=\h_i(\rr(\gamma_i))\\
  \dotgamma_i&=c_ie_i(\varphi)\,,\label{eq:10}
\end{align}
where~$e_i$ is the phase error in~\eqref{eq:5}.  The first thing we
notice about this system \eqref{eq:9}--\eqref{eq:10} is that it admits
a synchronized solution at any frequency $\omega^\star\!\in\n\real$
which belongs to the image of each function $\h_i\circ\rr$.  Indeed,
choose any initial states such that $\delta\varphi(0)=0$ and
$\h_i(\rr(\gamma_i(0)))=\omega^\star$ for each $i\n\in\n\ver$; then
$\dot{\varphi}_i\equiv\omega^\star$ and $\dotgamma_i\equiv 0$ for
every $i\n\in\n\ver$, which implies $\delta\varphi\equiv 0$.  However,
solutions starting from other initial conditions will not converge to
such synchronized trajectories, which means this system still fails to
meet the design goal of Definition~\ref{def:1}.

The second thing we notice about the system
\eqref{eq:9}--\eqref{eq:10} is that it is a Hamiltonian system with
position variables~$\varphi_i$ and momentum variables $\gamma_i/c_i$.
Indeed, assuming the phase coupling function~$f$ is odd, the
$n$-vector of phase errors~$e_i$ is the negative gradient of a
potential function of the phase vector~$\varphi$ (as we will show in
Section~\ref{sec:glob-lyap-funct}).  Furthermore, each function
$(\h_i\circ\rr)/c_i$, being a scalar function of a scalar variable,
is trivially the gradient of a potential function of its
argument~$\gamma_i$.  The sum of these potential functions is a
Hamiltonian associated with the dynamics \eqref{eq:9}--\eqref{eq:10}.
Moreover, as we will see in Section~\ref{sec:glob-lyap-funct}, if we
let the position variables be the phases~$\varphi_i$ measured relative
to an appropriate rotating frame, then the system is still
Hamiltonian, but now the Hamiltonian is proper and nonnegative and is
thus a Lyapunov function candidate.  If we can guarantee that
synchronized trajectories represent the only local minima of this
Lyapunov function, then we can perturb these Hamiltonian dynamics with
dissipation terms in the controller to achieve almost-global
synchronization.  This is the control design strategy we will pursue
in this paper.

\subsection{The PI loop filter: a perturbed Hamiltonian system}
\label{sec:loop-filter}

We propose two different perturbed Hamiltonian systems in this paper,
the first now and the second later on in
section~\ref{sec:dual-contr-anoth}.  For our first perturbation, we
add a proportional term to the integral control~\eqref{eq:8},
resulting in a proportional-integral (PI) loop filter of the form
\begin{align}
  \label{eq:12}
  \dotgamma_i&=c_i e_i(\varphi)\\
  u_i&=\rr(e_i(\varphi)+\gamma_i)\,.\label{eq:13}
\end{align}
Such a loop filter (without our nonlinear scaling function $\zeta$) can be found in
classical PLL designs~\cite{gar79} as well as in various coupled oscillator
systems~\cite{simspabarstr08,lun11}.  More generally, the error term
$e_i$ in~\eqref{eq:13} can be replaced with a scaled error term $\kk_i
e_i$, where the ``proportional gain'' $\kk_i$ is a positive constant
(or even a positive function of~$\gamma_i$), with virtually no change
in the convergence proof.  We have left this gain~$\kk_i$ out of the
analysis to simplify notation, but the flexibility it adds will be
important for tuning the performance of the system.  We might also
introduce a corresponding ``integral gain,'' but for the analysis this
can be absorbed into the constant~$c_i$.

To summarize, each controlled oscillator in the network has
second-order dynamics of the form
\begin{align}
  \label{eq:14}
  \dot{\varphi}_i&=\h_i(\rr(e_i(\varphi)+\gamma_i))\\
  \dotgamma_i&=c_i e_i(\varphi)\\
\intertext{with}
  \e_i(\varphi)&=\smashoperator{\sum_{(j,k)\in\N_i}}
   a_{k}f(\varphi_j-\varphi_i)\,.\label{eq:15}
\end{align}
The state space of each such oscillator is the cylinder
$\torus\!\times\n\real$.  We next present conditions under which this
oscillator system exhibits almost-global synchronization within a
given interval~$\I$.

\subsection{Assumptions and main result}
\label{sec:assumpt-main-result}

We assume that each vertex knows its neighbors in the graph~$\G$ (as
they need to exchange phase information), but otherwise the graph is
unknown.  However, we do assume the following:
\begin{enumerate}[itemsep=0.25em,label=\textbf{(A\arabic*)},%
                  align=left,leftmargin=*,series=assumptions]
\item the graph~$\G$ is connected, and there exists a known upper
  bound on the number~$n$ of vertices in~$\G$.\label{a:1}
\end{enumerate}
The frequency functions~$\h_i$ of the base oscillators are uncertain;
we merely assume that they satisfy the following:
\begin{enumerate}[assumptions]
\item the frequency functions $\h_i:\domg\to\real$ are all~$C^1$
  with positive derivatives $\h_i^\prime:\domg\to(0,\infty)$, and are
  such that\label{a:2}
  \begin{align}
    \label{eq:11}
    \bigcap_{i\in\ver}\h_i(\domg)&\neq\varnothing\,.
  \end{align}
\end{enumerate}
It is clear from the oscillator model~\eqref{eq:4} that~\eqref{eq:11}
is necessary for the existence of phase trajectories having
synchronized frequencies.  This condition~\eqref{eq:11} implies the
existence of some common interval of possible base oscillator
frequencies, and thus places an inherent limit on the extent to which
these oscillators can differ from each other.  Indeed, if one
oscillator can only produce frequencies in the kilohertz range and
another only in the megahertz range, then there is no possibility of
synchronization.

We next assume that our loop filter scaling function~$\rr$ satisfies:
\begin{enumerate}[assumptions]
\item the scaling function $\rr:\real\to\domg$ is~$C^1$ with
  positive derivative $\rr^\prime:\real\to(0,\infty)$, and
  is such that\label{a:3}
  \begin{align}
    \label{eq:16}
    \bigcap_{i\in\ver}\h_i(\rr(\real))&\neq\varnothing
    \intertext{and}
    \bigcup_{i\in\ver}\h_i(\rr(\real))&\subset\I\,,\label{eq:17}
  \end{align}
  where~$\I$ is the constraint interval from Definition~\ref{def:1}.
\end{enumerate}
It is clear from~\eqref{eq:4} that~\eqref{eq:17} constrains the
frequencies~$\dot{\varphi}_i$ to the interval~$\I$ as required by the
design goal in Definition~\ref{def:1}.  If~$\I$ is large enough to
contain the union of the images $\h_i(\domg)$, then we can always
satisfy assumption~\ref{a:3} by choosing the scaling function~$\rr$ to
be a diffeomorphism onto~$\domg$ (so that~\eqref{eq:11}
and~\eqref{eq:16} are the same).  However, if~$\I$ is not that large,
then assumption~\ref{a:3} states that we have found a solution to the
problem of designing~$\rr$ to satisfy both~\eqref{eq:16}
and~\eqref{eq:17} based on some \textit{a priori\/} knowledge about
the set of possible frequency functions~$\h_i$.  Such a design problem
could very well have no solution if~$\I$ is too small.

We let $f^\prime:\torus\to\real$ denote the derivative of the phase
coupling function~$f$ in the direction of the unit vector
field~$\unit$.  We make two assumptions on this function~$f$:
\begin{enumerate}[assumptions]
\item the phase coupling function $f:\torus\to\real$ is $C^1$
  and odd, that
  is, $f(-\theta)=-f(\theta)$ for all
  $\theta\n\in\n\torus$;\label{a:4}
\item there is a constant
  $b\in(0,\frac{\pi}{n-1}]$ such that $f^\prime\n(\theta)>0$ whenever
  $\cos(\theta)>\cos(b)$ and $f^\prime\n(\theta)<0$ whenever
   $\cos(\theta)<\cos(b)$,
  for any $\theta\n\in\n\torus$.
  \label{a:5}
\end{enumerate}
Assumptions \ref{a:4}--\ref{a:5} also appeared in
\cite{maltan10,maltan11} and play an important role in the convergence 
analysis: \ref{a:4} allows the interpretation of \eqref{eq:4} and \eqref{eq:7}-\eqref{eq:8}  
as a Hamiltonian system, while \ref{a:5} is needed to guarantee convergence
to the desired solution.

\begin{figure}[tp]
  \newcommand*{\fb}{0.5}
  \newcommand*{\fp}{2}
  \newcommand*{\fe}{0.53079}
  \newcommand*{\fP}{5}
  \newcommand*{\fE}{0.51490}
  \newcommand*{\Fp}{19}
  \newcommand*{\Fe}{0.50437}
  \centering
  \begin{tikzpicture}[trim axis left]
	\tikzstyle{every node}=[font=\footnotesize]

    \begin{axis}[  
    				scale only axis, 
				height=4cm,
				width=.5\columnwidth,
    				xmin=-3.5,xmax=3.5,
                 xtick={-3.14159,-\fb,0,\fb,3.14159},
                 xticklabels={$-\pi$,$-b$,$0$,$b$,$\pi$},
                 ytick={-0.4,-0.2,0,0.2,0.4},
                 legend pos=north west,
                 legend cell align=left,
                 grid=major,
                 major grid style={dotted}]
      \addplot[blue,thick,domain=-pi:pi,samples=100]
        {(1-cos(deg(\fb)))*sin(deg(x))/(1-cos(deg(\fb))*cos(deg(x)))};
      \addplot[green,thick,domain=-\fe:\fe,samples=100]
        {x*(1-(1/((\fp+1)*(\fb^\fp)))*abs(x)^\fp)};
      \addplot[green,thick,domain=\fe:pi,samples=2]
      {\fe*(1-(1/((\fp+1)*(\fb^\fp)))*\fe^\fp)+(1-(\fe/\fb)^\fp)*(x-\fe)};
      \addplot[green,thick,domain=-pi:-\fe,samples=2]
      {\fe*((1/((\fp+1)*(\fb^\fp)))*\fe^\fp-1)+(1-(\fe/\fb)^\fp)*(x+\fe)};
      \addplot[red,thick,domain=-\fE:\fE,samples=100]
        {x*(1-(1/((\fP+1)*(\fb^\fP)))*abs(x)^\fP)};
      \addplot[red,thick,domain=\fE:pi,samples=2]
      {\fE*(1-(1/((\fP+1)*(\fb^\fP)))*\fE^\fP)+(1-(\fE/\fb)^\fP)*(x-\fE)};
      \addplot[red,thick,domain=-pi:-\fE,samples=2]
      {\fE*((1/((\fP+1)*(\fb^\fP)))*\fE^\fP-1)+(1-(\fE/\fb)^\fP)*(x+\fE)};
      \addplot[cyan,thick,domain=-\Fe:\Fe,samples=100]
        {x*(1-(1/((\Fp+1)*(\fb^\Fp)))*abs(x)^\Fp)};
      \addplot[cyan,thick,domain=\Fe:pi,samples=2]
      {\Fe*(1-(1/((\Fp+1)*(\fb^\Fp)))*\Fe^\Fp)+(1-(\Fe/\fb)^\Fp)*(x-\Fe)};
      \addplot[cyan,thick,domain=-pi:-\Fe,samples=2]
      {\Fe*((1/((\Fp+1)*(\fb^\Fp)))*\Fe^\Fp-1)+(1-(\Fe/\fb)^\Fp)*(x+\Fe)};
      \legend{tanlock,degree $p=3$,,,degree $p=6$,,,degree $p=20$,,}
    \end{axis}
  \end{tikzpicture}
  \caption{Examples of phase-coupling functions $f$ when
    $b=0.5$.}
  \label{fig:2}
\end{figure}
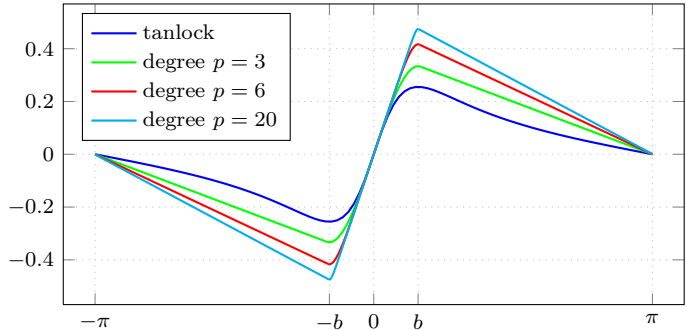

Note that to choose such a parameter $b$,
we must use our assumed knowledge in~\ref{a:1} of a known upper bound
on~$n$.  Examples of functions~$f$ which satisfy assumptions
\ref{a:4}--\ref{a:5} are shown in Fig.~\ref{fig:2} for $b\n=0.5$.  The
first example is given by the $C^\omega$ formula
\begin{align}
  f(\theta)
    &=\bigl[1-\cos(b)\bigr]
    \frac{\sin(\theta)}{1-\cos(b)\cos(\theta)}\,.
\end{align}
This function is related to the characteristic of certain ``tanlock''
phase comparators \cite{gar79}, and it generates the sinusoidal
coupling $f\n=\sin$ when $b\n=\frac{\pi}{2}$.  The other examples are
$C^1$ and piecewise polynomial of various degrees $p\geqslant 1$, each
having a derivative given by
$f^\prime\n(\theta)=1-(\abs{\theta}/b)^{p-1}$ on a certain arc
containing $[-b,b]$ and a constant derivative on its complement.  All
of these examples are normalized to have unit derivative at zero,
which means they should result in similar performance for small
deviations around a stable synchronized trajectory.  When~$b$ is small
(which we require when~$n$ is large), the magnitude of the derivative
of the tanlock function is small on the arc $[b,\pi]$ when
compared to the magnitude of the derivatives of the
piecewise-polynomial functions.  As a result, the piecewise-polynomial
functions might provide faster convergence to a synchronized state
when some initial phase differences are greater than~$b$
(due to their larger gains for large phase differences).

The final assumption is on the choice of the edge weights~$a_k$:
\begin{enumerate}[assumptions]
\item each edge weight~$a_k$ is chosen at random from a continuous
  probability distribution on the interval~$(0,\infty)$.
  \label{a:6}
\end{enumerate}
This assumption allows us to state that, with probability one, we
avoid an unknown zero-measure set of bad edge weight vectors
$a=[a_1\,\ldots\,a_m]^\tr\in\real^m$ for which our stability analysis
does not guarantee convergence.
\begin{theorem}\label{thm:1}
  Assume\/ {\em\ref{a:1}--\ref{a:6}}.  Then with probability one in
  the selection of edge weights in\/ {\em\ref{a:6}},
  the network of oscillators with vertex dynamics\/
  \eqref{eq:14}--\eqref{eq:15} achieves almost-global synchronization
  within\/~$\I$.
\end{theorem}

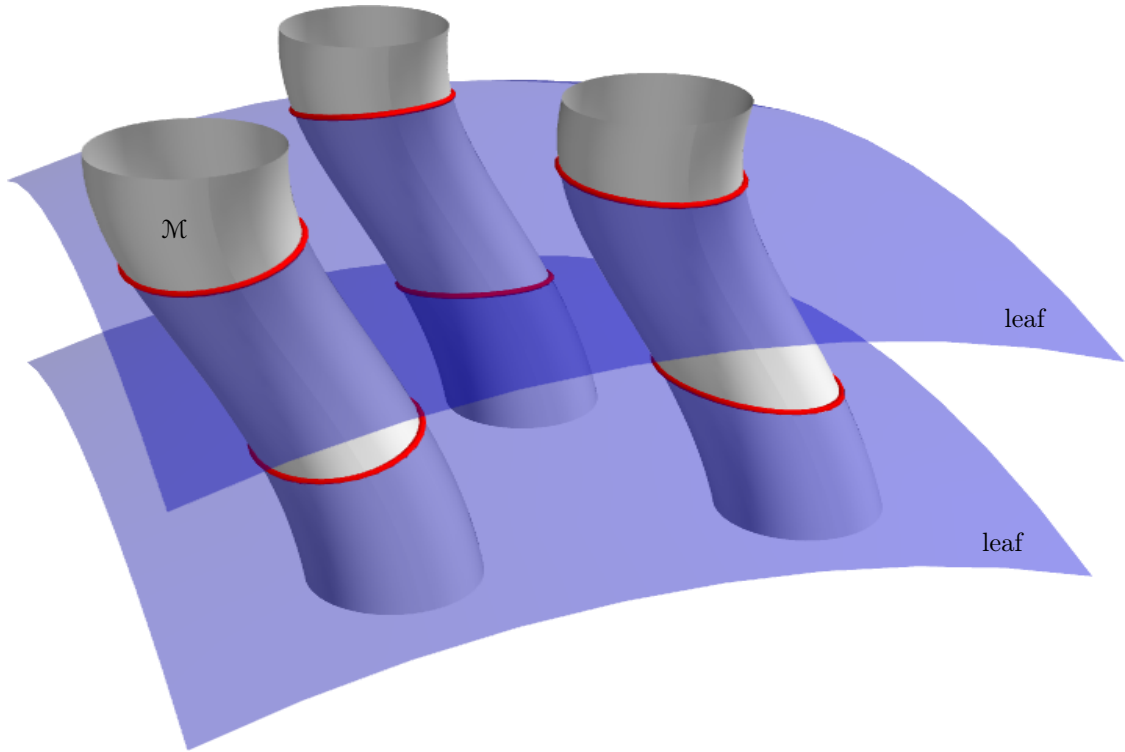
\begin{figure}
  \centering
  \begin{tikzpicture}
    \node[anchor=south west,inner sep=0] (image) at (0,0)
      {\includegraphics[width=0.9\columnwidth,height=.6\columnwidth]{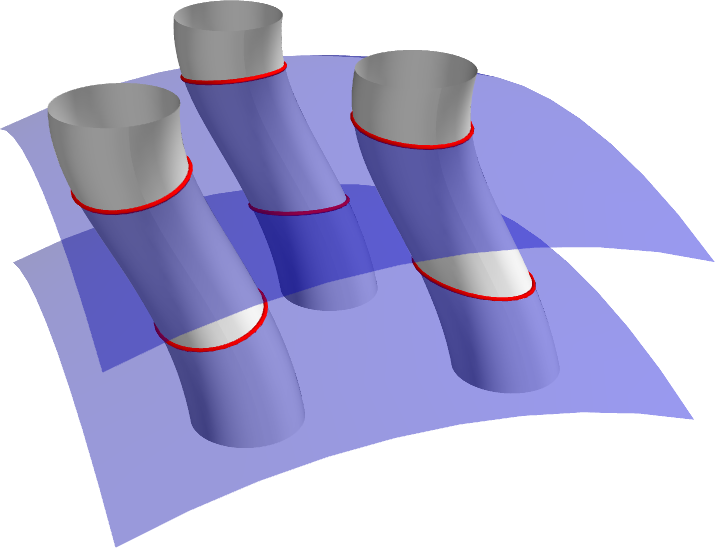}};
    \begin{scope}[x={(image.south east)},y={(image.north west)}]
      \node at (0.15,0.7) {$\M$};
      \node at (0.91,0.58) {leaf};
      \node at (0.89,0.28) {leaf};
    \end{scope}
  \end{tikzpicture}
  \caption{The invariant set $\M$ (shown in gray) is a collection of
    isolated cylinders.  The intersection of each invariant leaf of
    the foliation (shown in blue) with~$\M$ is a collection of
    isolated circles (shown in red).}
  \label{fig:3}
\end{figure}

Note that when $n\geqslant 4$, assumption~\ref{a:5} rules out the
sinusoidal phase-coupling function $f=\sin$.  In fact, the requirement
$b\n\leqslant\n\frac{\pi}{n-1}$ is within a factor of two of being
necessary.  Indeed, one can show that if $b>\frac{2\pi}{n}$, then for
some connected graphs having~$n$ vertices there exists a nonempty open
set of bad edge weight vectors such that trajectories of
\eqref{eq:14}--\eqref{eq:15} starting from a set of initial states
having positive measure achieve asymptotic frequency synchronization
but not asymptotic phase synchronization.

We will prove Theorem~\ref{thm:1} in Section~\ref{sec:proof}.  In this
proof, we first construct a Lyapunov function for the system, and then
apply the Krasovskii-LaSalle invariance theorem to show that all
trajectories converge to a certain invariant subset~$\M$ of the state
space $\torus^n\!\times\n\real^n$.  We then determine the structure of
this set~$\M$, showing first that it is the disjoint union of isolated
embeddings of the cylinder $\torus\!\times\n\real$.  We next show that
the state space admits a foliation into invariant $(2n-1)$-dimensional
submanifolds, and that the intersection of each leaf of this foliation
with the set~$\M$ is the disjoint union of isolated embeddings of the
circle~$\torus$ (as conceptualized in Fig.~\ref{fig:3}).  Each such
circle represents a periodic trajectory with a constant synchronized
frequency $\omega^\star\!\in\n\I$ (so that
$\dot{\varphi_i}\n\equiv\n\omega_i^\star$).  One of the cylinders of
the set~$\M$ is good, in that all of the periodic trajectories on this
cylinder are also synchronized in phase (so that $\delta\varphi\equiv
0$).  The remaining cylinders of~$\M$ (if any) are bad, in that all of
the periodic trajectories on these cylinders are out of phase.
Finally, using a local linearization analysis, we show that the
phase-synchronized trajectories are exponentially stable (relative to
each leaf of the foliation), whereas the out-of-phase trajectories are
exponentially unstable.  Because we show that these periodic
trajectories are \textit{isolated\/} on each leaf (a key step in the
proof), we conclude that almost all trajectories achieve asymptotic
synchronization in both frequency and phase.



\subsection{The dual controller: a different perturbed
Hamiltonian system}
\label{sec:dual-contr-anoth}

The approach we take in the controller design in
\eqref{eq:14}--\eqref{eq:15} is to perturb the Hamiltonian system
\eqref{eq:9}--\eqref{eq:10} by adding a dissipation term to the
dynamics of each position variable.  An alternative, dual approach is
to add a dissipation term to the dynamics of each momentum variable.
We will do this by including a second comparator in addition to the
phase comparator, that is, by introducing a second set of error
signals $z_i$ in addition to the phase errors~$e_i$ in~\eqref{eq:5}.
The resulting controller is more complex, but further analysis might
reveal that it has performance advantages over the controller
\eqref{eq:14}--\eqref{eq:15}.

For this dual controller, each pair of neighboring vertices agrees on
an orientation of the edge connecting them.  Thus each vertex~$i$ can
partition its neighbor set as $\N_i=\N_i^+\cup\N_i^-$, where
$(j,k)\in\N_i^+$ when~$i$ is the head of the edge $k=\{i,j\}$ and
$(j,k)\in\N_i^-$ when~$i$ is the tail of the edge~$k$.
Furthermore, we make an additional assumption on the frequency
functions~$\h_i$ of the base oscillators:
\begin{enumerate}[assumptions]
\item the interval~$\domg$ and the intervals $\h_i(\domg)$ for each
  $i\n\in\n\ver$ contain only positive real numbers.\label{a:7}
\end{enumerate}
In other words, we assume that the graph of each $\chi_i(\cdot)$ lies 
completely in the first quadrant, which is true for many analog oscillators.
This assumption guarantees that the functions $\rr$ and $\h_i$ take on
only positive values.  Consequently, for each ordered pair of
vertices~$ij$ we can define the ratio
\begin{align}
  \label{eq:20}
  \rho_{ij}(\gamma_i,\gamma_j)&=
  \frac{\rr(\gamma_j)\h_i(\rr(\gamma_i))}%
  {\rr(\gamma_i)\h_j(\rr(\gamma_j))}>0\,.
\end{align}
Using~\eqref{eq:9}, that is, using the base oscillator
model~\eqref{eq:4} together with the integrator loop filter
\eqref{eq:7}--\eqref{eq:8}, we can calculate~\eqref{eq:20} as
\begin{align}
  \label{eq:21}
  \rho_{ij}(\gamma_i,\gamma_j)
  &=\frac{\rr(\gamma_j)}{\rr(\gamma_i)}
  \Bigl[\frac{d\varphi_j}{d\varphi_i}\Bigr]^\inv\,.
\end{align}
We then define a new error signal~$z_i$ for vertex~$i$ as
\begin{align}
  z_i(\gamma)&=\smashoperator{\sum_{(j,k)\in\N_i^+}}
          d_k\bigl[\rr(\gamma_j)-\rho_{ij}(\gamma_i,\gamma_j)
          \rr(\gamma_i)\bigr]\notag\\
     &\qquad
         +\smashoperator{\sum_{(j,k)\in\N_i^-}}
          d_k\bigl[\rho_{ij}^\inv(\gamma_i,\gamma_j)
          \rr(\gamma_j)-\rr(\gamma_i)\bigr]\,,
  \label{eq:22}
\end{align}
where $\gamma=[\gamma_1\,\ldots\,\gamma_n]^\tr\n\in\n\real^n$ is the
vector of loop filter states and the constants $d_{k}$ are an
additional set of positive edge weights.  For vertex~$i$ to compute
this error variable~$z_i$, it must receive the signal
$\rr(\gamma_j)$ from each of its neighbors.  Also, we see
from~\eqref{eq:21} that it must calculate the derivative of each
neighbor's phase~$\varphi_j$ with respect to its own
phase~$\varphi_i$, a calculation which does not require knowledge of
the global time variable~$t$.

Using these new error signals~$z_i$, we define the dual perturbed
Hamiltonian system as
\begin{align}
  \label{eq:23}
  \dot{\varphi}_i&=\h_i(\rr(\gamma_i))\\
  \dotgamma_i&=c_ie_i(\varphi)+c_iz_i(\gamma)\,.\label{eq:24}
\end{align}
Thus we have retained the simple integrator loop filter in
\eqref{eq:7}--\eqref{eq:8}, but have changed its input from~$e_i$ to
the sum $e_i+z_i$.  This amounts to modifying Fig.~\ref{fig:1} by
putting a ``$\gamma$-comparator'' in parallel with the phase
comparator.  

Although the at first sight the physical interpretation of $z_i$
is obscured by the implementation details, substituting \eqref{eq:20} and
\eqref{eq:23} into \eqref{eq:22} gives
\begin{align}
  z_i(\gamma)&=\smashoperator{\sum_{(j,k)\in\N_i}}\eta_k(\gamma)
  \bigl[\chi_j(\zeta(\gamma_j))-\chi_i(\zeta(\gamma_i))\bigr]\,\nonumber\\
  &= \smashoperator{\sum_{(j,k)\in\N_i}}\eta_k(\gamma)
  \bigl[\dot\varphi_j-\dot\varphi_i\bigr]\,  \label{eq:67}
\end{align}
with 
\begin{align}
  \label{eq:66}
  \eta_k(\gamma)&=
  d_k\dfrac{\rr(\gamma_j)}{\chi_j(\zeta(\gamma_j))}\,,\text{ where } j
  \text{ is the tail of edge }k,
\end{align}
for each $k\n\in\n\edg$.  Therefore \eqref{eq:22} can be interpreted as an 
indirect  estimate of the frequency mismatch \eqref{eq:67} between oscillator
$i$ and its neighbors.
As we will see in the proof of the following theorem in
Section~\ref{sec:proof-theor-refthm:1} below, this frequency error terms $z_i$ 
are the right perturbation that need to be added to the dynamics of the momentum 
variables to guarantee the dissipation of the Hamiltonian energy.
\begin{theorem}\label{thm:3}
  Assume {\em\ref{a:1}--\ref{a:7}}.  Then with probability one in the
  selection of edge weights in\/ {\em\ref{a:6}}, the network of
  oscillators with vertex dynamics\/ \eqref{eq:22}--\eqref{eq:24}
  achieves almost-global synchronization within\/~$\I$.
\end{theorem}
The proof of this theorem is similar to that of Theorem~\ref{thm:1},
and we will highlight the main differences in
Section~\ref{sec:proof-theor-refthm:1}.

\section{Proof of Theorem~2}
\label{sec:proof}

For any $p\n\in\n\nat$, we let $\0_p$ and $\1_{\n p}$ denote the
column vectors of~$p$ zeros and~$p$ ones, respectively, and we
let~$I_p$ denote the $p\n\times\n p$ identity matrix.  Any
$\ell\n\times\n p$ matrix with integer elements will represent either
an $\real$-linear map from~$\real^p$ to~$\real^\ell$ or a
$\intgr$-linear map from~$\torus^p$ to~$\torus^\ell$, depending on
context.  In particular, if~$M$ is an $\ell\n\times\n p$ matrix with
integer elements and $t\mapsto x(t)$ is a curve in~$\torus^p$, then
$t\mapsto y(t)=Mx(t)$ is a curve in~$\torus^\ell$ (with~$M$
representing a $\intgr$-linear map), and furthermore
$\doty(t)=M\dotx(t)$ for all~$t$ (with~$M$ now representing an
$\real$-linear map).

We first introduce some notation for writing the overall coupled
dynamics \eqref{eq:14}--\eqref{eq:15} in a compact form.  We define
state vectors
\begin{align}
  \varphi&=\bigl[\varphi_1\,\ldots\,\varphi_n\bigr]^\tr&
  \in\;&\torus^n\\
  \gamma&=\bigl[\gamma_1\,\ldots\,\gamma_n\bigr]^\tr&\in\;&\real^n
\end{align}
along with the following diagonal matrices:
\begin{align}
  \label{eq:25}
  A&=\diag\{a_{1},\,\dotsc,\,a_{m}\}
  &\in\;&\real^{m\times m}\\
  C&=\diag\{c_1,\dotsc,c_n\}&\in\;&\real^{n\times n}\,.
\end{align}
For convenience we define $\g_i=\h_i\circ\rr$ for each~$i$, and we
note from assumptions \ref{a:2}--\ref{a:3} that each~$\sigma_i$
is~$C^1$ with a positive derivative.  We next define mappings
$F:\torus^m\to\real^m$ and $\gv:\real^n\to\real^n$ by
\begin{align}
  F(\theta)&=
  \begin{bmatrix}
    f(\theta_1)&
    \ldots&
    f(\theta_m)
  \end{bmatrix}^\tr&\in\;&\real^m\\
  \gv(y)&=
  \begin{bmatrix}
    \g_1(y_1)&
    \ldots&
    \g_n(y_n)
  \end{bmatrix}^\tr&\in\;&\real^n\,,\label{eq:26}
\end{align}
where $\theta=[\theta_1\,\ldots\,\theta_m]^\tr\n\in\n\torus^m$ and
$y=[y_1\,\ldots\,y_n]^\tr\n\in\n\real^n$.  Note that~$\gv$ is a $C^1$
diffeomorphism onto its image $\gv(\real^n)$.  Finally, we let
$B\n\in\n\{-1,0,1\}^{n\times m}$ be an oriented incidence matrix for
the graph~$\G$.  Using the fact from~\ref{a:4} that~$f$ is odd, we can
write the dynamics \eqref{eq:14}--\eqref{eq:15} as
\begin{align}
  \dot{\varphi}&=\gv(-B\n AF(B^\tr\!\varphi)
   +\gamma)\label{eq:27}\\
  \dotgamma&=-C\n B\n AF(B^\tr\!\varphi)\,.
  \label{eq:28}
\end{align}
For convenience we define $\w:\torus^n\to\real^n$ by
\begin{align}
  \label{eq:29}
  \w(\varphi)&=\bigl[e_1(\varphi)\,\ldots\,e_n(\varphi)\bigr]^\tr
  =\;-B\n AF(B^\tr\!\varphi)\,,
\end{align}
which enables us to write
\eqref{eq:27}--\eqref{eq:28} as
\begin{align}
  \dot{\varphi}&=\gv(\w(\varphi)
   +\gamma)\label{eq:30}\\
  \dotgamma&=C\w(\varphi)\,.
  \label{eq:31}
\end{align}
Note that $\1_n^\tr\w\equiv 0$ because $\1_n^\tr B=0$ (a property of
any oriented incidence matrix).  This property also implies
$\w(\varphi)=\w(\varphi+\1_n\theta)$ for any $\varphi\n\in\n\torus^n$
and any $\theta\n\in\n\torus$.

\subsection{Global Lyapunov function}
\label{sec:glob-lyap-funct}

In this section we construct a Lyapunov function for the system
\eqref{eq:30}--\eqref{eq:31}, which is the Hamiltonian we described in
Section~\ref{sec:towards-synchr-hamil} for the system
\eqref{eq:9}--\eqref{eq:10}.

Because $f$ is odd from~\ref{a:4}, the integral of the 1-form
$f\n\cdot\n\ip{\unit,\cdot}$ around any smooth closed curve
in~$\torus$ is zero.  Thus this 1-form is the differential of a smooth
function $\Psi:\torus\to\real$, which is unique up to an additive
constant (which we choose so that the minimum value of~$\Psi$
on~$\torus$ is zero).  Therefore $\frac{d}{dt}(\Psi\circ x)\equiv
f(x)\dotx$ for any curve $x:\real\to\torus$.
We then define $V:\torus^m\to[0,\infty)$ as the sum
\begin{align}
  \label{eq:32}
  V(\theta)&=\sum_{k\in\edg}a_k\Psi(\theta_k)\,,
\end{align}
where
$\theta=[\theta_1\,\ldots\,\theta_m]^\tr\n\in\n\torus^m$.
It follows that
\begin{align}
  \label{eq:33}
  \frac{d}{dt}V(B^\tr\!\varphi)
  &=F^\tr\n(B^\tr\!\varphi)AB^\tr\n\dot{\varphi}
   =-\w^\tr\n(\varphi)\dot{\varphi}\,.
\end{align}
From~\eqref{eq:16} there exists a frequency $\gstar\n\in\bigcap_i
\g_i(\real)$.  Therefore the function $W:\real^n\to[0,\infty)$
defined as
\begin{align}
  \label{eq:34}
  W(\gamma)&=\sum_{i\in\ver}\frac{1}{c_i}
  \int_{\g_i^\inv(\gstar)}^{\gamma_i}
  \bigl[\g_i(s)-\gstar\bigr]ds
\end{align}
is proper (because~$\g_i$ is increasing) and has a derivative given by
\begin{align}
  \label{eq:35}
  \frac{d}{dt}
  W(\gamma)&=\bigl[\gv^\tr\n(\gamma)
    -\gstar\1_n^{\n\tr}\bigr]C^\inv\dotgamma\,.
\end{align}
We now obtain a Lyapunov function by adding~\eqref{eq:32}
and~\eqref{eq:34}: we define $U:\torus^n\!\times\real^n\to[0,\infty)$
by
\begin{align}
  \label{eq:36}
  U(\varphi,\gamma)&=V(B^\tr\!\varphi)+W(\gamma)\,,
\end{align}
which is a proper function on $\torus^n\!\times\real^n$.  Because
$\1_n^\tr\w\equiv 0$, the derivative of~\eqref{eq:36} along
trajectories of \eqref{eq:30}--\eqref{eq:31} is
\begin{align}
  \dotU&=-\w^\tr\n(\varphi)\gv(\w(\varphi)
   +\gamma)+\gv^\tr\n(\gamma)\w(\varphi)\notag\\
  &=-\w^\tr\n(\varphi)\n\int_0^1\!\!
  \gv^\prime\n
  (\w(\varphi)s+\gamma)\,ds\cdot
  \w(\varphi)\,,
  \label{eq:37}
\end{align}
where $\gv^\prime\n(\cdot)$ denotes the diagonal Jacobian matrix
\begin{align}
  \gv^\prime\n(\gamma)
  &=\diag\bigl\{\g_1^\prime(\gamma_1),\,
  \dotsc,\,\g_n^\prime(\gamma_n)\bigr\}\,.\label{eq:38}
\end{align}
Because $\gv^\prime\n(\cdot)$ is positive definite, we have
$\dotU\leqslant 0$ and $\dotU=0$ if and only if $\w(\varphi)=0$.  It
follows from the Krasovskii-LaSalle invariance theorem that all
trajectories of the system \eqref{eq:30}--\eqref{eq:31} converge to
the largest invariant set~$\M$ contained within the set
$\Phi\n\times\real^n$, where $\Phi\subset\torus^n$ denotes the set
\begin{align}
  \label{eq:39}
  \Phi&=\w^\inv(\{0\})=\bigl\{\varphi\in\torus^n\,
             :\,B\n AF(B^\tr\!\varphi)=0\bigr\}\,.
\end{align}
Note that because $B^\tr\n\1_n=0$, this set~$\Phi$ has the symmetry
property $\Phi=\Phi+\1_n\n\torus$.  Also, it follows from~\ref{a:4} that
$F(0)=0$, which means $\Phi$ contains all points of the form
$\1_n\theta$ for $\theta\n\in\n\torus$.  The next step in the proof is
to investigate the structure of this set~$\M$.

\subsection{Structure of the set\/~$\M$}
\label{sec:structure-set-m}

In this section we show that the largest invariant set~$\M$ contained
within the set $\Phi\n\times\real^n$ is $\M=\Phi\n\times\n\Gamma$,
where~$\Phi$ in~\eqref{eq:39} is the zero set of~$\w$ and
$\Gamma\subset\real^n$ is the set
\begin{align}
  \label{eq:40}
  \Gamma&=\bigl\{\gamma\in\real^n\,
             :\,B^\tr\n \gv(\gamma)=0\bigr\}\,.
\end{align}
We begin by exploring the structure of these sets~$\Phi$ and~$\Gamma$.
First, we show that $\Phi$ is the disjoint union of isolated
embeddings of the circle~$\torus$.  Next, we show that
$\Gamma\n=\alpha(\real)$, where $\alpha:\real\to\real^n$ is a $C^1$
curve in~$\real^n$.  Moreover, if we let $q\in\real^n$ denote the unit
vector in the direction of $C^\inv\1_n$, then this curve~$\alpha$ is
such that $q^\tr\n\alpha(\cdot)$ is the identity map on~$\real$.
Finally, we show that $\M=\Phi\n\times\n\Gamma$, that is, that $\M$ is
the disjoint union of isolated embeddings of the cylinder, as
illustrated in Fig.~\ref{fig:3}.

We partition~$\torus^n$ using the following matrices:
\begin{align}
  \label{eq:41}
  R&=
  \begin{bmatrix}
    \0^\tr_{n-1}\\I_{n-1}
  \end{bmatrix}&\in\;&\{0,1\}^{n\times(n-1)}\\
  S&=
  \begin{bmatrix}
    -\1_{n-1}^\tr\\I_{n-1}
  \end{bmatrix}&\in\;&\{-1,0,1\}^{n\times(n-1)}\,.\label{eq:42}
\end{align}
Clearly the state~$\varphi$ satisfies the identity
\begin{align}
  \label{eq:43}
  \varphi&=\1_n\varphi_1+RS^\tr\!\varphi\,,
\end{align}
which defines the direct sum $\torus^n\n=\1_n\n\torus\oplus
R\torus^{n-1}$.  Here the first summand represents the first component
angle and the second one represents the remaining angles measured
relative to the first.  Note that $S^\tr\!R=I_{n-1}$, that
$S^\tr\n\1_n\n=0$, and that $B^\tr\!=\n B^\tr\!RS^\tr$.  Also,
because~$\G$ is connected from~\ref{a:1}, we have $\rank(B)=n-1$ and
thus the columns of $B^\tr\! R$ are independent.

Because $\varphi\n\in\n\Phi$ if and only if
$RS^\tr\!\varphi\n\in\n\Phi$, it follows from~\eqref{eq:43} that
$\Phi=RS^\tr\n\Phi+\1_n\n\torus$.  It also follows that for any
$\mu\in\torus^{n-1}$, we have $\mu\n\in\n S^\tr\n\Phi$ if and only if
$R\mu\n\in\n\Phi$, that is, if and only if $\w(R\mu)=0$.  We next show
that the points in the set $S^\tr\n\Phi$ are isolated, which implies
that $\Phi$ is the disjoint union of isolated embeddings of the
circle~$\torus$.

Using the above partition of~$\torus^n$, we define
two symmetric matrix functions $L:\torus^{n-1}\to\real^{n\times n}$ and
$L^\flat:\torus^{n-1}\to\real^{(n-1)\times(n-1)}$ by
\begin{align}
  \label{eq:44}
  L(\mu)&=B\n AF^\prime\n(B^\tr\!R\mu)B^\tr&\text{and}&&
  L^\flat(\mu)&=R^\tr\n L(\mu) R
\end{align}
for $\mu\in\torus^{n-1}$, where $F^\prime\n(\cdot)$ denotes the
diagonal Jacobian matrix
\begin{align}
  F^\prime\n(\theta)
  &=\diag\bigl\{f^\prime\n(\theta_1),\,
    \dotsc,\,f^\prime\n(\theta_m)\bigr\}
\end{align}
with $\theta=[\theta_1\,\ldots\,\theta_m]^\tr\in\torus^m$.
Here~$L(\mu)$ represents a weighted Laplacian matrix for the
graph~$\G$ in which the weights can have positive, negative, or zero
values.  Also, because $B=SR^\tr\n B$ we have
\begin{align}
  \label{eq:45}
  L(\mu)R&=SL^\flat(\mu)&&\text{and}&L(\mu)&=SL^\flat(\mu)S^\tr
\end{align}
for all $\mu\n\in\n\torus^{n-1}$.  Note that $L(\mu)$ is congruent to
the block diagonal matrix $\diag\{0,L^\flat(\mu)\}$.  The proof of the
following theorem is in Appendix~\ref{sec:proof-theor-refthm:2}:
\begin{theorem}\label{thm:4}
  There is a set\/ $\mathscr{Z}\subset\real^m$ having zero Lebesgue
  measure such that if\/ $a=[a_1\,\ldots\,a_m]^\tr\not\in\mathscr{Z}$,
  then the matrix\/ $L^\flat(\mu)$ in\/ {\em\eqref{eq:44}} is
  invertible for all\/ $\mu\in\n S^\tr\n\Phi$.
\end{theorem}
\begin{corollary}\label{cor:1}
  If\/ $a\not\in\mathscr{Z}$ then the points in\/
  $S^\tr\n\Phi$ are isolated.
\end{corollary}
\begin{IEEEproof}
  Define the mapping $P:\torus^{n-1}\n\to\real^{n-1}$ by setting
  $P(\mu)=R^\tr\! B\n AF(B^\tr\! R\mu)$ so that
  $P^\inv(\{0\})=S^\tr\n\Phi$.  The Jacobian matrix for~$P$ is just
  $L^\flat(\mu)$, which by Theorem~\ref{thm:4} is invertible for all
  $\mu\n\in\n S^\tr\n\Phi$.  The result follows from the inverse
  function theorem.
\end{IEEEproof}

We do not provide a method for computing the set~$\mathscr{Z}$ of bad
edge weight vectors; instead, we rely on the random edge weight
selection in~\ref{a:6} to avoid this zero-measure set.  It is clear
from~\eqref{eq:39} and~\eqref{eq:44} that this set~$\mathscr{Z}$
depends only on the graph~$\G$ (through~$B$) and the phase coupling
function~$f$ (through~$F$).  In principle we could remove the
dependence on the unknown graph by taking the countable union of all
such sets~$\mathscr{Z}$ over all connected graphs to be our
zero-measure set of bad edge weights, but this larger set would still
be difficult to compute.  Thus from now on we will assume
$a\not\in\mathscr{Z}$, which happens with probability one according
to~\ref{a:6}

To explore the structure of the set $\Gamma$ in~\eqref{eq:40}, we
partition~$\real^n$ by defining $q\in\real^n$ and
$Q\in\real^{n\times(n-1)}$ as
\begin{align}
  q&=\frac{C^\inv\1_n}{\norm{C^\inv\1_n}}&\text{and}&&
  Q&=CS\n\bigl(S^\tr\n C^2S\bigr)^{\n-\frac{1}{2}}\,.\label{eq:46}
\end{align}
Using the fact that $\1_n^\tr S=0$, it is straightforward to show that
the $n\n\times\n n$ matrix $[q\;\:Q]$ is orthogonal; thus the
state~$\gamma$ satisfies the identity
\begin{align}
  \label{eq:47}
  \gamma&=qq^\tr\n\gamma+QQ^\tr\n\gamma\,,
\end{align}
which defines the direct sum $\real^n=q\real\oplus Q\real^{n-1}$ via
orthogonal projections.  We define the open interval $J\subset\real$
as
\begin{align}
  \label{eq:48}
  J&=\bigcap_{i\in\ver}\g_i(\real)\,,
\end{align}
which is nonempty from~\eqref{eq:16}.  Because $s\1_n\in\gv(\real^n)$
for all $s\in J$, we can define the function $\beta:J\to\real$ as
\begin{align}
  \label{eq:49}
  \beta(s)&=\sum_{i\in\ver}q_i\g_i^\inv(s)
           =q^\tr\gv^\inv\n(s\1_n)\,,
\end{align}
where the constants~$q_i$ are the components of the vector~$q$.  Note
that each term in the sum in~\eqref{eq:49} is strictly increasing
in~$s$.  Suppose $\{s_\ell\}$ is a sequence in~$J$ such that
$s_\ell\to\sup J$ as $\ell\to\infty$.  Then there exists
$i\n\in\n\ver$ such that $s_\ell\to\sup\g_i(\real)$ as
$\ell\to\infty$, which means $\g_i^\inv(s_\ell)\n\to\n\infty$ as
$\ell\n\to\n\infty$.  Thus $\beta(s_\ell)\n\to\n\infty$ as
$\ell\n\to\n\infty$, which means $\sup\beta(J)=\infty$.  Similar
reasoning yields $\inf\beta(J)=-\infty$, and we conclude that~$\beta$
is invertible.  Therefore we can define a $C^1$ curve
$\alpha:\real\to\real^n$ as follows:
\begin{align}
  \label{eq:50}
  \alpha(r)&=\gv^\inv(\beta^\inv(r)\1_n)\,,
\end{align}
which satisfies $q^\tr\n\alpha(r)=\beta(\beta^\inv(r))=r$ for all
$r\n\in\n\real$.  Because $\rank(B)=n-1$, it follows
from~\eqref{eq:40} that $\gamma\n\in\n\Gamma$ if and only if
$\gv(\gamma)=x\1_n$ for some $x\in\n J$, that is, if and only if
\begin{align}
  \label{eq:51}
  \gamma&=\gv^\inv(x\1_n)=\gv^\inv(\beta^\inv(\beta(x))\1_n)
         =\alpha(\beta(x))
\end{align}
for some $x\in\n J$.  Therefore
$\Gamma\n=\alpha(\beta(J))=\alpha(\real)$.

Our remaining task in this section is to show that
$\M=\Phi\n\times\n\Gamma$.  First we calculate the derivative
of~$\w(\varphi)$ in~\eqref{eq:29} using~\eqref{eq:30}:
\begin{align}
  \label{eq:52}
  \dot{\w}&=-B\n AF^\prime\n(B^\tr\!\varphi)B^\tr\gv(\w(\varphi)
   +\gamma)\notag\\&=-L(S^\tr\!\varphi)\,\gv(\w(\varphi)
   +\gamma)\,.
\end{align}
We see from \eqref{eq:30}--\eqref{eq:31} that
$\dot{\varphi}=\gv(\gamma)$ and $\dotgamma=0$ on the set
$\Phi\n\times\real^n$, which means the second derivative $\ddotvarphi$
is zero:
\begin{align}
  0=\ddotvarphi&=-\gv^\prime\n(\gamma)
  L(S^\tr\!\varphi)\,\gv(\gamma)\notag\\
  &=-\gv^\prime\n(\gamma)
  SL^\flat(S^\tr\!\varphi)S^\tr\n\gv(\gamma)
  \label{eq:53}
\end{align}
on $\Phi\n\times\real^n$.  Now $\gv^\prime\n(\cdot)$ in~\eqref{eq:38}
is positive definite, $L^\flat(S^\tr\!\varphi)$ is invertible
from Theorem~\ref{thm:4}, and~$S$ in~\eqref{eq:42} has independent
columns; therefore~\eqref{eq:53} implies
$S^\tr\n\gv(\gamma)=0$. Because $B^\tr\!=\n B^\tr\!RS^\tr$, this in
turn implies $B^\tr\n\gv(\gamma)\n=0$, and we conclude that
$\M\subset\Phi\n\times\n\Gamma$.  We observe that
$\Phi\n\times\n\Gamma$ is itself invariant under the dynamics
\eqref{eq:30}--\eqref{eq:31}, and it follows that
$\M=\Phi\n\times\n\Gamma$.

To summarize the results of this section, we have found that
$(\varphi,\gamma)\n\in\n\M$ if and only if $\w(\varphi)\n=\n 0$ and
$\gamma\n=\n\alpha(r)$ for some $r\n\in\n\real$, in which case
$r=q^\tr\n\gamma$.  Furthermore, the zero set of~$\w$ in~\eqref{eq:39}
is $RS^\tr\n\Phi+\1_n\n\torus$, where the points in $S^\tr\n\Phi$ are
isolated.

\subsection{Global analysis: frequency synchronization}
\label{sec:global-analysis}

To study the synchronization properties of our system
\eqref{eq:30}--\eqref{eq:31}, we first observe that
$q^\tr\n\dotgamma\equiv 0$, which means the state space admits a
foliation whose leaves are the invariant manifolds
$\torus^n\!\times\Xi_r$, where
\begin{align}
  \label{eq:54}
  \Xi_r&=qr+Q\real^{n-1}=
  \bigl\{\gamma\in\real^n\,:\,q^\tr\n\gamma=r\bigr\}
\end{align}
for a constant parameter $r\in\real$.  Note that because
$\Gamma\n=\alpha(\real)$ and $q^\tr\n\alpha(\cdot)$ is the identity
map, the intersection $\Gamma\cap\Xi_r$ is the singleton
$\{\alpha(r)\}$.  Thus the intersection of each invariant leaf
$\torus^n\!\times\Xi_r$ with $\M=\Phi\n\times\n\Gamma$ is just
$\Phi\times\{\alpha(r)\}$, which we have shown to be the disjoint
union of isolated embeddings of the circle~$\torus$ (depicted as red
circles in Fig.~\ref{fig:3}).

We now fix $r\in\real$ and examine the dynamics on
$\torus^n\!\times\Xi_r$, noting that $\gamma\equiv qr+QQ^\tr\n\gamma$
on this invariant manifold.  Now that~$r$ is fixed, we will write
$\alpha_r=\alpha(r)$.  Because $\alpha_r=qr+QQ^\tr\n\alpha_r$, we have
\begin{align}
  \label{eq:55}
  \gamma-\alpha_r&\equiv QQ^\tr\n(\gamma-\alpha_r)
\end{align}
on $\torus^n\!\times\Xi_r$.
We next define the projected state variables
\begin{align}
  \label{eq:56}
  w_1&=S^\tr\!\varphi\in\torus^{n-1}&\text{and}&&
  w_2&=Q^\tr\n(\gamma-\alpha_r)\in\real^{n-1}\,.
\end{align}
Taking time derivatives of these variables, and using~\eqref{eq:55}
and the fact that $\w(\varphi)=\w(Rw_1)$, we obtain
\begin{align}
  \label{eq:57}
  \dot{w}_1&=S^\tr\n\gv(\w(Rw_1)+Qw_2+\alpha_r)\\
  \dot{w}_2&=Q^\tr\! C\w(Rw_1)\,.\label{eq:58}
\end{align}
This is an autonomous system in the projected states $(w_1,w_2)$, and
its equilibria are precisely all points of the form $(\mu^\star,0)$
for vectors $\mu^\star\n\in\n S^\tr\n\Phi$.  Each equilibrium
represents a frequency-synchronized solution of
\eqref{eq:30}--\eqref{eq:31} with
$\dot{\varphi}\n\equiv\n\beta^\inv(r)\1_n$ and
$\gamma\n\equiv\n\alpha_r$.  The equilibrium with $\mu^\star\n=0$
represents a phase-synchronized trajectory with $\delta\varphi\equiv
0$, and all other equilibria represent out-of-phase trajectories.

Furthermore, each trajectory of the system
\eqref{eq:57}--\eqref{eq:58} converges to an equilibrium, which means
each trajectory of the system \eqref{eq:30}--\eqref{eq:31} achieves
asymptotic frequency synchronization.  Indeed, we have shown that all
trajectories of the system \eqref{eq:30}--\eqref{eq:31} converge to
the set $\M=\Phi\n\times\n\Gamma$ in forward time.  Therefore $\gamma$
converges to the singleton set $\Gamma\cap\Xi_r=\{\alpha_r\}$, which
means $w_2$ converges to zero.  In addition, $\varphi$ converges to
the set~$\Phi$, which means $w_1$ converges to the set $S^\tr\n\Phi$.
Corollary~\ref{cor:1} states that points in $S^\tr\n\Phi$ are
isolated, and we conclude that $w_1$ converges to one of these points.

In the next section, we will perform a local linearization analysis at
each equilibrium of the system \eqref{eq:57}--\eqref{eq:58} to
determine its stability.

\subsection{Local analysis: phase synchronization}
\label{sec:local-analysis}

We compute the linear approximation of the dynamics
\eqref{eq:57}--\eqref{eq:58} at an equilibrium $(\mu^\star,0)$ as
follows:
\begin{align}
  \label{eq:59}
  \dot{w}_1&\approx -S^\tr\n\gv^\prime\n(\alpha_r) L(\mu^\star)
            R(w_1-\mu^\star)+S^\tr\n\gv^\prime\n(\alpha_r)Qw_2\\
  \dot{w}_2&\approx -Q^\tr\! CL(\mu^\star)R(w_1-\mu^\star)\,,
\end{align}
or using~\eqref{eq:45},
\begin{align}
  \label{eq:60}
  \dot{w}_1&\approx -S^\tr\n\gv^\prime\n(\alpha_r)
            SL^\flat(\mu^\star)
            (w_1-\mu^\star)+S^\tr\n\gv^\prime\n(\alpha_r)Qw_2\\
  \dot{w}_2&\approx -Q^\tr\! CSL^\flat(\mu^\star)(w_1-\mu^\star)\,.
\end{align}
If we define the $(n-1)\n\times\n(n-1)$ matrices
\begin{align}
  X&=S^\tr\n\gv^\prime\n(\alpha_r) S>0\\
  Y&=Q^\tr\! C S=\bigl(S^\tr\n C^2
     S\bigr)^{\n\frac{1}{2}}>0\label{eq:61}\\
  Z&=S^\tr\n\gv^\prime\n(\alpha_r)Q\,,\label{eq:62}
\end{align}
then we can write this approximation more compactly as
\begin{align}
  \dot{w}_1&\approx -X\n
  L^\flat(\mu^\star)(w_1-\mu^\star)+Zw_2\label{eq:63}\\ 
  \dot{w}_2&\approx -Y\n
  L^\flat(\mu^\star)(w_1-\mu^\star)\,.\label{eq:64}
\end{align}
Note from~\eqref{eq:46} and~\eqref{eq:61} that $Q=CSY^\inv$ so that
$S=C^\inv QY$, which means~$Z$ in~\eqref{eq:62} satisfies
\begin{align}
  \label{eq:65}
  Y^\inv Z&=Q^\tr\! C^\inv \gv^\prime\n(\alpha_r)Q>0\,.
\end{align}
Thus the linearization \eqref{eq:63}--\eqref{eq:64} satisfies the
assumptions in the following theorem, whose proof is in
Appendix~\ref{sec:proof-theor-refthm:3}:
\begin{theorem}\label{thm:5}
  Let\/ $\Lambda\in\real^{2p\times 2p}$ be the block matrix
  \begin{align}
    \Lambda&=
    \begin{bmatrix*}[r]
      -X\n L & Z \\ -Y\n L & 0\,
    \end{bmatrix*},
  \end{align}
  where\/ $L,X,Y,Z\in\real^{p\times p}$ satisfy:
  \begin{enumerate}
  \item $L$ is symmetric,
  \item $X\n+\n X^\tr\n>0$,
  \item $Y$ is symmetric and invertible, and
  \item $\displaystyle Y^\inv\n Z$ is symmetric with\/
    $\displaystyle Y^\inv\n Z\geqslant 0$.
  \end{enumerate}
  If\/ $L$ has a negative eigenvalue, then\/ $\Lambda$ has
  an eigenvalue with a positive real part.  If instead\/ $Z$
  is invertible and\/ $L>0$, then\/ $\Lambda$ is Hurwitz.
\end{theorem}
We can complete the local linearization analysis by applying
Theorem~\ref{thm:5} together with the following theorem:
\begin{theorem}\label{thm:6}
  The Laplacian\/ $L(\mu)$ has\/ $n-1$ positive eigenvalues for\/
  $\mu\n=\n0$, and it has at least one negative eigenvalue for any
  nonzero\/ $\mu\n\in\n S^\tr\n\Phi$.
\end{theorem}
The proof of Theorem~\ref{thm:6}, which can be found in
\cite{mallada2013synchronization}, relies on the properties
of the phase coupling function~$f$ defined in Assumption \ref{a:5}.
The basic intuition is that for 
$\mu^\star\not=0$, choosing $b$ small enough guarantees a 
negative eigenvalue in  $L(\mu^\star)$ independently of the choice 
of $\mu^\star\in S^T\Phi$.

Let us now consider an equilibrium $(\mu^\star,0)$ of the nonlinear
system \eqref{eq:57}--\eqref{eq:58}.  Recall that $L(\mu^\star)$ is
congruent to $\diag\{0,L^\flat(\mu^\star)\}$; thus Theorem~\ref{thm:6}
together with Sylvester's law of inertia imply that
$L^\flat(\mu^\star)\n>\n0$ when $\mu^\star\n=\n0$ and that
$L^\flat(\mu^\star)$ has a negative eigenvalue for all nonzero
$\mu^\star\n\in\n S^\tr\n\Phi$.  Therefore if $\mu^\star\n=0$, which
represents an in-phase steady-state solution, then it follows from
Theorem~\ref{thm:5} that this equilibrium is exponentially stable.
Likewise, if $\mu^\star\n\neq 0$, which represents an out-of-phase
steady-state solution, then this equilibrium is exponentially
unstable.  Because all out-of-phase equilibria are both isolated and
exponentially unstable, we conclude (say
from~\cite[Proposition~1]{fre13}, for example) that the set
$\E_r\subset\torus^n\!\times\Xi_r$ of initial states from which
trajectories converge to out-of-phase steady-state solutions has zero
measure in~$\torus^n\!\times\Xi_r$ (regarded here as a
$(2n-1)$-dimensional manifold).  It then follows from Tonelli's
theorem that the set $\E=\bigcup_{r\in\real}\E_r$ has zero measure in
$\torus^n\!\times\n\real^{n}$.  In other words, the system achieves
asymptotic phase synchronization from almost every initial state.

\section{Proof of Theorem~3}
\label{sec:proof-theor-refthm:1}

We define the diagonal matrix function
\begin{align}
  \label{eq:68}
  H(\gamma)&=\diag\bigl\{\eta_{1}(\gamma),\,\dotsc,
  \,\eta_{m}(\gamma)\bigr\}
  &\in\;&\real^{m\times m}\,,
\end{align}
so that we can write the dynamics \eqref{eq:23}--\eqref{eq:24} as
\begin{align}
  \label{eq:69}
  \dot{\varphi}&=\gv(\gamma)\\
  \dotgamma&=C\w(\varphi)-C\n BH(\gamma)B^\tr\gv(\gamma)\,.
  \label{eq:70}
\end{align}
Because $\1_n^\tr\w\equiv 0=\1_n^\tr B$, the derivative
of~\eqref{eq:36} along trajectories of \eqref{eq:69}--\eqref{eq:70} is
\begin{align}
  \label{eq:71}
   \dotU&=-\gv^\tr\n(\gamma)BH(\gamma)B^\tr\gv(\gamma)\,.
\end{align}
Now $H(\cdot)$ is positive definite; hence $\dotU\leqslant 0$ and
$\dotU=0$ if and only if $B^\tr\gv(\gamma)=0$.  It follows from the
Krasovskii-LaSalle invariance theorem that all trajectories of the
system \eqref{eq:69}--\eqref{eq:70} converge to the largest invariant
set~$\M$ contained within the set $\torus^n\!\times\n\Gamma$, where
$\Gamma\subset\real^n$ is from~\eqref{eq:40}.  We next show that again
$\M=\Phi\n\times\n\Gamma$, just as in the proof of
Theorem~\ref{thm:1}.

It follows from~\eqref{eq:70} that this systems also satisfies
$q^\tr\n\dotgamma\equiv 0$ and thus admits a foliation whose leaves
are the invariant manifolds $\torus^n\!\times\Xi_r$.  Because
$B^\tr\gv(\gamma)$ is the constant zero on the set~$\Gamma$, its
derivative is zero on~$\M$:
\begin{align}
  0&\equiv B^\tr\gv^\prime\n(\gamma)\dotgamma\equiv
     B^\tr\gv^\prime\n(\gamma)QQ^\tr\n\dotgamma\notag\\
   &\equiv B^\tr\!RS^\tr\gv^\prime\n(\gamma)CS\n
     \bigl(S^\tr\n C^2S\bigr)^{\n-\frac{1}{2}}Q^\tr\n\dotgamma\,.
  \label{eq:72}
\end{align}
Because the columns of $B^\tr\! R$ and~$S$ are independent and because
the diagonal matrix $\gv^\prime\n(\gamma)C$ is positive definite, it
follows that $Q^\tr\n\dotgamma\equiv 0$ and thus $\dotgamma\equiv 0$
on~$\M$.  It then follows from~\eqref{eq:70} that $\w(\varphi)\equiv
0$ on~$\M$, and we conclude that $\M\subset\Phi\n\times\n\Gamma$.  We
observe that $\Phi\n\times\n\Gamma$ is itself invariant under the
dynamics \eqref{eq:69}--\eqref{eq:70}, and it follows that
$\M=\Phi\n\times\n\Gamma$.

We now fix $r\in\real$ and examine the dynamics on the invariant
manifold $\torus^n\!\times\Xi_r$.  Using~\eqref{eq:55}, we see that
the derivatives of the projected state variables~\eqref{eq:56} along
trajectories of \eqref{eq:69}--\eqref{eq:70} are
\begin{align}
  \label{eq:73}
  \dot{w}_1&=S^\tr\n\gv(Qw_2+\alpha_r)\\
  \dot{w}_2&=Q^\tr\! C\w(Rw_1)\notag\\
    &\qquad-Q^\tr\! C\n BH(Qw_2+\alpha_r)B^\tr\gv(Qw_2+\alpha_r)\,.
  \label{eq:74}
\end{align}
As in the proof of Theorem~\ref{thm:1}, the equilibria of this system
are precisely all points of the form $(\mu^\star,0)$ for vectors
$\mu^\star\n\in\n S^\tr\n\Phi$.  The equilibrium with $\mu^\star\n=0$
represents a fully synchronized trajectory with $\delta\varphi\equiv
0$, and all other equilibria represent out-of-phase,
frequency-synchronized trajectories.  Because
$\M=\Phi\n\times\n\Gamma$, because $\Gamma\cap\Xi_r$ is a singleton,
and because the points of $S^\tr\n\Phi$ are isolated, we see that each
trajectory of the system \eqref{eq:73}--\eqref{eq:74} converges to an
equilibrium.


We now perform the local linearization analysis at each
equilibrium $(\mu^\star,0)$ of the system
\eqref{eq:73}--\eqref{eq:74}. Similarly, to \eqref{eq:63}--\eqref{eq:64} we can 
approximate \eqref{eq:73}--\eqref{eq:74} around $(\mu^\star,0)$ using 
\begin{align}
  \dot{w}_1&\approx Zw_2\label{eq:84}\\
  \dot{w}_2&\approx -Y\!L^\flat(\mu^\star)(w_1-\mu^\star)-Y\!L^\flat_\eta Zw_2\label{eq:85}
\end{align}
where
\begin{align}
L^\flat_\eta(\alpha_r)&=R^T L_\eta(\alpha_r) R &\text{ and }&& L_\eta(\alpha_r)&=BH(\alpha_r)B^T,
\end{align}
with the matrices  $L_\eta(\alpha_r)$ and $L^\flat_\eta(\alpha_r)$ also satisfying a condition analogous to \eqref{eq:45}.
  

We can now use a dual version of Theorem~\ref{thm:5} based on \cite[Theorem~5, Corollary~2]{carsch62} 
which we state below.
\begin{lemma}\label{lem:337}
  Suppose\/ $\Lambda\in\cmplx^{q\times q}$ has no eigenvalues on the
  imaginary axis, suppose\/ $H\!\in\cmplx^{q\times q}$ is an
  invertible Hermitian matrix, and suppose\/ $\Lambda
  H+H\!\Lambda^{\!\star}\geqslant 0$, where\/ $\Lambda^{\!\star}$
  denotes the conjugate transpose of\/ $\Lambda$.  Then the number of
  eigenvalues of\/ $\Lambda$ having positive real part is the same as
  the number of positive eigenvalues of\/~$H$.
\end{lemma}

The dual version of Theorem \ref{thm:5} is then:
\begin{theorem}\label{thm:7}
  Let\/ $\Lambda\in\real^{2p\times 2p}$ be the block matrix
  \begin{align}
    \Lambda&=
    \begin{bmatrix*}[r]
      0\; & Z \\ -Y\!L_1 & -Y\!L_2 Z
    \end{bmatrix*},\label{eq:58}
  \end{align}
  where\/ $L_1,L_2,Y,Z\in\real^{p\times p}$ satisfy:
  \begin{enumerate}
  \item $L_1$ is symmetric and invertible,
  \item $L_2$ is symmetric with $L_2>0$,
  \item $Y$ is symmetric and invertible, and
  \item $\displaystyle Y^\inv\!Z$ is symmetric with\/
    $\displaystyle Y^\inv\!Z>0$.
  \end{enumerate}
  If\/ $L_1$ has a strictly negative eigenvalue, then\/ $\Lambda$ has
  an eigenvalue with a strictly positive real part.  If instead\/
  $L_1>0$, then\/ $\Lambda$ is Hurwitz.
\end{theorem}
The proof is in Appendix \ref{sec:proof-theor-refthm:7}.

Thus, we can use Theorem~\ref{thm:6} together with Theorem~\ref{thm:7} to show that, for the system \eqref{eq:84}--\eqref{eq:85},  the equilibrium with $\mu^\star=0$ is exponentially stable and any other equilibria with $\mu^\star\not=0$ is exponentially unstable. The rest of the proof follows that of Theorem~\ref{thm:1}.

\section{Concluding remarks}
\label{sec:concluding-remarks}

We have presented two distributed controllers for the phase and
frequency synchronization of heterogeneous nonlinear oscillators, and
we have shown that each guarantees almost-global convergence on
arbitrary connected graphs. 
Our solutions can be readily implemented using analog oscillators and 
phase comparators, and our analysis holds under very general 
assumptions on the system. In particular, unlike most existing work, 
we neither require the set of admissible frequencies to be unbounded, nor
assume any special network topology.

\appendix

\subsection{Proof of Theorem~4}
\label{sec:proof-theor-refthm:2}

  We let $\mathscr{T}$ denote the finite collection of all
  $m\n\times\n m$ diagonal matrices
  $\Delta=\diag\{\delta_1,\dotsc,\delta_m\}$ such that for all
  $k\in\edg$, either $\delta_k=f^\prime\n(0)$ or
  $\delta_k=f^\prime\n(\pi)$.  For each such matrix
  $\Delta\in\mathscr{T}$, we define the closed set
\begin{align}
  \mathscr{P}_{\n\Delta}&=\bigl\{a\in\real^m\,:\,\det\bigl(R^\tr\!
  B\diag(a)\Delta B^\tr\! R\bigr)=0\bigr\}\,,
\end{align}
where $\diag(a)=\diag\{a_1,\dotsc,a_m\}$ denotes the diagonal matrix
whose diagonal entries are the~$m$ elements of~$a$.  Now~$\Delta$ is
invertible by assumption~\ref{a:5}, and furthermore the columns of
$B^\tr\! R$ are independent; it follows that
$\mathscr{P}_{\n\Delta}\neq\real^m$ (take $\diag(a)=\Delta^{\n\inv}$),
which means~$\mathscr{P}_{\n\Delta}$ is a closed algebraic set having
zero measure.  Thus
\begin{align}
  \mathscr{P}&=\bigcup_{\Delta\in\mathscr{T}}\mathscr{P}_{\n\Delta}
\end{align}
is also a closed algebraic set having zero measure.  Therefore the set
$\mathscr{O}=\real^m\setminus\mathscr{P}$ is a nonempty open
semialgebraic set.  Next we define the mapping
$\mathscr{H}:\torus^{n-1}\times\mathscr{O}\to\real^{n-1}$ by
\begin{align}
  \mathscr{H}(\mu,a)&=R^\tr\! B\diag(a)F(B^\tr\! R\mu)\,.
  \label{eq:75}
\end{align}
The Jacobian matrix of~$\mathscr{H}$ is
\begin{align}
  D\mathscr{H}(\mu,a)&=
  \begin{bmatrix}
    \dfrac{\partial\mathscr{H}}{\partial \mu}(\mu,a)
    &\dfrac{\partial\mathscr{H}}{\partial a}(\mu,a)
  \end{bmatrix}\,,
\end{align}
where
\begin{align}
  \dfrac{\partial\mathscr{H}}{\partial \mu}(\mu,a)&=
    R^\tr\! B\diag(a)F^\prime\n(B^\tr\! R\mu)B^\tr\! R
    \label{eq:76}\\[0.5em]
  \dfrac{\partial\mathscr{H}}{\partial a}(\mu,a)&=
    R^\tr\! B\diag\bigl(F(B^\tr\! R\mu)\bigr)\,.
\end{align}
If we define the matrix
\begin{align}
  \mathscr{J}(\mu,a)&=
  \begin{bmatrix}
    I\\
    \!\!\!\!\diag(a)\diag\bigl(F(B^\tr\! R\mu)\bigr)^{\n+}\quad\\
    \quad\cdot\bigl[F^\prime\n(0)
  -F^\prime\n(B^\tr\! R\mu)\bigr]B^\tr\! R
  \end{bmatrix},
\end{align}
where $(\cdot)^+$ denotes the Moore–Penrose pseudoinverse, then
\begin{align}
  D\mathscr{H}(\mu,a)\cdot\mathscr{J}(\mu,a)&=
  R^\tr\! B\diag(a)\Delta(\mu)B^\tr\! R\,,\label{eq:77}
\end{align}
where~$\Delta(\mu)$ is the diagonal matrix
\begin{align}
  \Delta(\mu)&=F^\prime\n(B^\tr\! R\mu)
  +\diag\bigl(F(B^\tr\! R\mu)\bigr)
  \diag\bigl(F(B^\tr\! R\mu)\bigr)^{\n+}\notag\\
  &\qquad\qquad\cdot\bigl[F^\prime\n(0)
  -F^\prime\n(B^\tr\! R\mu)\bigr]\,.\label{eq:78}
\end{align}
It follows from assumptions \ref{a:4}--\ref{a:5} that $f(\theta)=0$ if
and only if $\theta\n\in\n\{0,\pi\}$, so for any
$\mu\n\in\n\torus^{n-1}$, the matrix $\Delta(\mu)$ in~\eqref{eq:78}
belongs to $\mathscr{T}$.  It follows from the definition
of~$\mathscr{O}$ that the matrix in~\eqref{eq:77} is invertible, and
we conclude that $D\mathscr{H}(\mu,a)$ has rank $n-1$ for all
$(\mu,a)\in\torus^{n-1}\n\times\mathscr{O}$.  Thus
$\mathscr{H}\pitchfork\{0\}$,%
\footnote{If~$M$ and~$N$ are smooth manifolds, if $f:N\to M$ is
  smooth, and if $S$ is an embedded submanifold of $M$, then $f$ is
  \textit{transverse to\/} $S$, written $f\pitchfork S$, when for
  every $p\in f^\inv(S)$ we have $T_{f(p)}M=T_{f(p)}S + df_p(T_pN)$,
  where $df_p$ denotes the differential of~$f$ at~$p$.}  and it
follows from the parametric transversality theorem\footnote{This
  transversality theorem is based on Sard's theorem, which holds only
  for sufficiently smooth functions.  One can verify that in our case,
  this application of Sard's theorem is valid when the phase coupling 
  function~$f$ is merely~$C^1$.}  \cite[Theorem~6.35]{lee12} that there exists a
set $\mathscr{Y}\subset\mathscr{O}$ having zero measure such that if
\mbox{$a\in\mathscr{O}\setminus\mathscr{Y}$} then
$\mathscr{H}_a\pitchfork\{0\}$, where 
$\mathscr{H}_a$ denotes the mapping
\mbox{$\mu\mapsto\mathscr{H}(\mu,a)$}.  Choose
$\mathscr{Z}=\mathscr{P}\cup\mathscr{Y}$; we have thus shown that for
all $a\in\real^m\setminus\mathscr{Z}$, if~$\mu$ is such that
$\mathscr{H}(\mu,a)=0$, then the matrix in~\eqref{eq:76} is
invertible.

Suppose $a\in\real^m$ is the edge weight vector so that $A=\diag(a)$,
suppose $a\n\in\n\real^m\n\setminus\mathscr{Z}$, and suppose
$\mu\n\in\n S^\tr\n\Phi$.  Then $R\mu\n\in\n\Phi$ which means
$\mathscr{H}(\mu,a)=0$, and it follows that $L^\flat(\mu)$
in~\eqref{eq:44}, which is the same as the matrix in~\eqref{eq:76}, is
invertible.

\subsection{Proof of Theorem~6}
\label{sec:proof-theor-refthm:3}

  Because both $Y$ and $\displaystyle Y^\inv\n Z$ are symmetric
  matrices, we see that the following holds:
  \begin{align}
    Y^\inv\n ZY&=\bigl(Y^\inv\n Z\bigr)^{\n\tr}Y=Z^\tr Y^\inv Y=Z^\tr\,.
    \label{eq:79}
  \end{align}
  We will investigate the stability of the zero solution of the system
  $\dot{x}=\Lambda x$.  Partition the state as
  $x=[x_1^\tr\;\;x_2^\tr]^\tr$ with $x_1,x_2\in\real^p$, and consider
  the quadratic function
  \begin{align}
    \Upsilon(x)&=\Upsilon(x_1,x_2)=x_1^\tr\n Lx_1+x_2^\tr Y^\inv\n Zx_2
  \end{align}
  whose derivative along trajectories of $\dotx=\Lambda x$ is
  \begin{align}
    \dot{\Upsilon}(x)\n&=-2x_1^\tr\n L\n X\n Lx_1\n+2x_1^\tr\n LZx_2
             -2x_2^\tr Y^\inv\n ZY\n Lx_1\notag\\
    &=-x_1^\tr\n L(X\n+\n X^\tr)Lx_1\leqslant 0\,.\label{eq:80}
  \end{align}
  First suppose~$L$ has an eigenvalue $\lambda<0$ with unit
  eigenvector~$v$.  The standard Chetaev instability conditions
  \cite[Theorem~4.3]{kha02} do not guarantee \textit{exponential\/}
  instability, so we will instead exploit the linear-quadratic
  structure to strengthen the Chetaev result.  The cone
  \begin{align}
    \mathscr{C}&=\bigl\{x\in\real^{2p}\,:\,\Upsilon(x)<0\bigr\}\,,
  \end{align}
  is nonempty because
  \begin{align}
    \Upsilon(v,0)&=v^\tr\n Lv
    =\lambda v^\tr\n v=\lambda<0\,.\label{eq:81}
  \end{align}
  For any $x\in\mathscr{C}$ we have
  \begin{align}
    0&\leqslant x_2^\tr Y^\inv\n Zx_2<-x_1^\tr\n Lx_1
    =\abs[\big]{x_1^\tr\n Lx_1}\,,
  \end{align}
  which means
  \begin{align}
    \abs{\Upsilon(x)}&\leqslant\abs[\big]{x_1^\tr\n Lx_1}
    +x_2^\tr Y^\inv\n Zx_2<2\abs[\big]{x_1^\tr\n Lx_1}\,.
  \end{align}
  Let $L^+$ denote the Moore–Penrose pseudoinverse of~$L$, which
  is nonzero because~$L$ has a nonzero eigenvalue~$\lambda$.  Then
  we have
  \begin{align}
    \abs[\big]{x_1^\tr\n Lx_1}&=\abs[\big]{x_1^\tr\n LL^+\n Lx_1}
    \leqslant\norm{L^+}\cdot\norm{Lx_1}^2\,.
  \end{align}
  Let $\sigma>0$ denote the minimum eigenvalue of
  $X\n+\n X^\tr$, so that
  \begin{align}
    \dot{\Upsilon}(x)&=-x_1^\tr\n L(X\n+\n X^\tr)Lx_1
    \leqslant-\sigma\norm{Lx_1}^2\,.
  \end{align}
  It follows that for any $x\in\mathscr{C}$ we have
  \begin{align}
    \dot{\Upsilon}(x)&<\frac{\sigma}{2\norm{\displaystyle L^+}}
    \Upsilon(x)<0\,.
  \end{align}
  Therefore~$\mathscr{C}$ is positively invariant and $\Upsilon(x(t))$
  exhibits exponential growth towards $-\infty$ as $t\to\infty$ from
  any initial state $x(0)\in\mathscr{C}$.
  Hence~$\norm{x(t)}=\norm[\big]{\exp(\Lambda t)x(0)}$ exhibits
  exponential growth as $t\to\infty$ when $x(0)\in\mathscr{C}$, and we
  conclude that $\Lambda$ must have an eigenvalue with a positive real
  part.

  Next suppose that $Z$ is invertible and $L>0$; then the largest
  invariant set contained in the set $\{\dot{\Upsilon}(x)=0\}$ is the
  origin $x=0$.  In this case~$\Upsilon$ is positive definite, hence
  asymptotic stability follows from the Krasovskii-LaSalle invariance
  theorem.

 \subsection{Proof of Theorem~9}
 \label{sec:proof-theor-refthm:7}

For the instability part of Theorem~\ref{thm:5}, we do not require the
matrices~$L_1$ or~$Z$ to be invertible, even though they happen to be
so in our application.  However, the above Chetaev-style proof of
Theorem~\ref{thm:5} does not carry over to a proof of the instability
part of Theorem~\ref{thm:7}; instead, we will apply
Lemma~\ref{lem:337}.
To apply this lemma to matrix~$\Lambda$ in~\eqref{eq:58}, we must show
that~$\Lambda$ has no eigenvalues on the imaginary axis.  Suppose
$\Lambda v=\ii \lambda v$ for some $\lambda\in\real$ and
$v\in\cmplx^{2p}$, and partition~$v$ as $v=[v_1^\tr\;\;v_2^\tr]^\tr$
with $v_1,v_2\in\cmplx^p$.  Then
\begin{align}
  Zv_2&=\ii\lambda v_1\label{eq:3378}\\
  -Y\!L_1v_1-Y\!L_2 Zv_2&=\ii\lambda v_2\,.\label{eq:3379}
\end{align}
If $\lambda=0$, then $v_1=v_2=0$ because $L_1$, $Y$, and~$Z$ are all
invertible.  Otherwise we solve for~$v_2$ in~\eqref{eq:3378} and
substitute into~\eqref{eq:3379} to obtain
\begin{align}
  Y\bigl[\lambda^2Y^\inv\!Z^\inv\!-L_1-\ii\lambda L_2\bigr]v_1&=0\,.
  \label{eq:3380}
\end{align}
If we write $v_1=x+\ii y$ for $x,y\in\real^p$, then
because~$Y$ is invertible we obtain
\begin{align}
  0&=
  \bigl[\lambda^2Y^\inv\!Z^\inv\!-L_1-\ii\lambda L_2\bigr]
  (x+\ii y)\notag\\
  &=\bigl[Xx+\lambda L_2y\bigr]
   +\ii\bigl[-\lambda L_2x+Xy\bigr]\,,
  \label{eq:3381}
\end{align}
where we have defined $X=\lambda^2Y^\inv\!Z^\inv\!-L_1$.  Setting the
imaginary part to zero in~\eqref{eq:3381} and solving for~$x$ yields
\begin{align}
  x&=\tfrac{1}{\lambda}L_2^\inv\!Xy\,.
  \label{eq:3382}
\end{align}
Setting the real part to zero in~\eqref{eq:3381} and
substituting~\eqref{eq:3382} for~$x$ yields
\begin{align}
  0&=\bigl[\tfrac{1}{\lambda}XL_2^\inv\!X
    +\lambda L_2\bigr]y\,.
  \label{eq:3383}
\end{align}
Now~\eqref{eq:3375} implies $ZY=Y\!Z^\tr$, which means
\begin{align}
  (Y^\inv\!Z^\inv)^\tr&=\bigl[(ZY)^\inv\bigr]^\tr
  =\bigl[(ZY)^\tr\bigr]^\inv\notag\\
  &=\bigl[Y\!Z^\tr\bigr]^\inv=\bigl[ZY\bigr]^\inv\notag\\
  &=Y^\inv\!Z^\inv\,,
\end{align}
and we conclude that~$X$ is symmetric.  Thus the matrix
$XL_2^\inv\!X+\lambda^2 L_2$ is symmetric and positive definite, and
therefore~\eqref{eq:3383} implies $y=0$.  Hence $x=0$
from~\eqref{eq:3382}, which means $v_1=0$ and thus also $v_2=0$
from~\eqref{eq:3378}.  Therefore, if $\Lambda v=\ii \lambda v$ for some
$\lambda\in\real$ and $v\in\cmplx^{2p}$ then $v=0$, and we conclude
that~$\Lambda$ has no eigenvalues on the imaginary axis.

Let $H\in\real^{2p\times 2p}$ be the block diagonal matrix
\begin{align}
  H&=
  \begin{bmatrix*}[c]
    -L_1^\inv & 0 \\
    0 & -Z^\inv Y
  \end{bmatrix*}\,,
\end{align}
which is symmetric and invertible by assumption.  Then a
straightforward calculation yields
\begin{align}
  \Lambda H+H\!\Lambda^{\!\tr}&=
  \begin{bmatrix*}[c]
    0 & 0 \\
    0 & 2Y\!L_2 Y
  \end{bmatrix*}\geqslant 0\,,
\end{align}
and the result follows from Lemma~\ref{lem:337}.
\printbibliography

\end{document}